\documentstyle[12pt,amsfonts]{article}
\def \k{\C [[\hbar ]]}    
\def \g{\tilde{\frak g}}
\def \dbar{\overline{\partial}}

\def \C{{\mathbb C}}

\def \R{{\mathbb R}}

\def \A{{\mathbb A}}
\def \fA{\A^{\hbar}}
\def \vfA{\hat{\A}^{\hbar}}
\def \W{{\mathbb W}}
\def \fx{\hat{x}}
\def \feta{\hat{\eta}}
\def \fz{\hat{z}}

\def \fxi{\hat{\xi }}

\newtheorem{th}{Theorem}[section]

\newtheorem{definition}[th]{Definition}

\newtheorem{proposition}[th]{Proposition}

\newtheorem{lemma}[th]{Lemma}

\newtheorem{corollary}[th]{Corollary}

\newtheorem{remark}[th]{Remark}

\newtheorem{example}[th]{Example}
\newtheorem{assumption}[th]{Assumption}
\newtheorem{notation}[th]{Notation}

\def \App{{\mathbb A}^{\hbar}}

\title{Deformations of symplectic Lie algebroids, deformations of holomorphic symplectic structures, and index theorems}

\author{ $Ryszard \ Nest^1$ and $Boris \ Tsygan^2$
 }

\date{}

\begin{document}

\maketitle

\tableofcontents

%

%



\section{Introduction}

A deformation quantization of a smooth manifold $M$ is by definition a formal multiplication law on the space $C^{\infty}(M)[[\hbar]]$
$$
f*g = fg + \sum_{k\geq 1} (i\hbar )^k P_k (f,g),
$$   
where $*$ is an associative $\hbar$-linear product satisfying 
$$
\frac{1}{i\hbar } (f*g -g*f) = \{ f,g \} +O(\hbar ).
$$
One usually requires $P_k$ to be local, i.e. bidifferential, expressions in $f$ 
and $g$. It is also convenient to assume that
$$
1*f =f*1 =f .
$$
\cite{BFFLS}.

A powerful theorem of Kontsevich \cite{Konts} states that there is a bijection between the set of isomorphism classes of deformations of $C^{\infty}(M)$ and the set of equivalence classes of formal Poisson structures (i.e., of formal series $\varpi = \sum _{k \geq 0} \hbar ^k \varpi ^k$ satisfying $[\varpi,\varpi]=0$ where $\{, \;\}$ is the Schouten-Nijenhuis bracket). In this paper, we classify deformations for a class of Poisson structures for which explicit methods of Fedosov work. Using those methods, we prove an index theorem for this class of deformations.

In the case when $M$ is a symplectic manifold it is known that deformation 
quantizations always exist \cite{LW} and are classified by the 
points of the affine space 
$
\frac{1}{i \hbar}\omega +H^2 (M, \C [[\hbar ]]) 
$
 \cite{NTfamilies}, \cite{De 2}, \cite{LW}. 
A simple geometric construction of deformation quantization of a symplectic 
manifold was given in \cite{F} by Fedosov.

 Fedosov's methods are well suited for 
study of a more general class of Poisson manifolds. In this paper we 
generalize both them and the classification results of \cite{NTfamilies} to the 
case of Poisson structures associated with symplectic Lie 
algebroids (this generality was suggested to us by A. Weinstein who also
independently carried out the construction of Fedosov connections in \cite{A.W}).

\medskip

A Lie algebroid (\cite{pradines}, \cite{?}, \cite{bb}) is a vector bundle $E$ over $M$ whose sheaf of sections is a sheaf of Lie algebras
and a morphism of bundles 
$$
\rho : E\rightarrow TM 
$$
satisfying
$$
\rho [\xi ,\eta ] =[\rho (\xi ),\rho (\eta )]
$$
and 
$$
[f\xi ,\eta ] = f [\xi , \eta ]- L_{\rho (\eta )} (f) \xi
$$
for $f$ in $C^{\infty} (M)$ and $\xi ,\eta $ in $\Gamma (M,E)$.

Given a Lie algebroid, one can define $E$-differential forms
$$
{}^E \Omega^* (M) =\Gamma (M, \Lambda E^* )
$$ 
and the de Rham differential
$$
d: {}^E \Omega^* (M) \rightarrow {}^E \Omega^{*+1} (M)
$$
(using the classical Cartan formula). A symplectic Lie algebroid is a Lie 
algebroid $(E,[ \ ,\ ] ,\rho )$ together with a non-degenerate closed 
element of ${}^E \Omega^2 (M)$. One can define an $E$-connection on a vector bundle $F$ as an operator $\Gamma (M,F) \rightarrow {}^E\Omega^1(M, F)$ satisfying standard properties. The typical examples of symplectic Lie algebroids are as follows
\begin{itemize}
\item {\it Symplectic manifolds} \\
The Lie algebroid in question is just the sheaf of vector fields on $M$ and
the deformation theory is determined by the second cohomology group $H^2 (M,
\C [[\hbar ]])$.
\item {\it Constant rank Poisson structures} \\
The Lie algebroid is given by the sheaf of vector fields tangent to the
leaves of a Poisson foliation and
the deformation theory is determined by the second cohomology group $H^2 (M,
\underline{O}_{hor} [[\hbar ]])$ (with coefficients in the sheaf of leafwise
constant functions). 
\item {\it Symplectic structures with logarithmic singularities along  submanifolds of codimension
    one} \\
The Lie algebroid is given by the sheaf of vector fields whose restriction to the submanifold is tangent to this submanifold, and
the deformation theory is determined by the second cohomology group ${}^b H^2
(M, \C [[\hbar ]])$, the de Rham cohomology of differential forms with
logarithmic singularities.
\item {\it Manifolds with corners}\\
The Lie algebroid is given by a subsheaf of the sheaf of vector fields whose 
restrictions 
to the boundary are tangent to the strata at the boundary. The 
deformation theory is determined by the second cohomology group ${}^b H^2
(M, \C [[\hbar ]]$, the de Rham cohomology of differential forms with logarithmic
singularities along the boundary.
\item {\it Compactified cotangent bundles}\\
For a manifold $X$, one can compactify its cotangent bundle $T^*X$ by adding the cosphere bundle $S^*X$ at infinity to get the closed ball bundle ${\overline{B}}^*X$. Let $E$ be the Lie algebroid of fields on ${\overline{B}}^*X$ whose restriction to $S^*X$ is tangent to the fibres $S^*_x$, $x \in X$. One can show that the standard symplectic form on $T^*X$ extends to an $E$-symplectic form on ${\overline{B}}^*X$.
\item {\it Complex symplectic manifolds}\\
The Lie algebroid is given by the sheaf of vector fields of type (1,0), and
the deformation theory involves both deformation of the holomorphic structure
and deformation of the product. In this case, there are natural obstructions to
construction of deformations. These obstructions are related to the Hodge spectral sequence. When
they vanish, the deformations are again classified by their characteristic
class, but the space of characteristic classes allowed is a proper (non
affine) submanifold of $H^2_{DR} (M)[[\hbar]]$.
\end{itemize}

\medskip

In the context of a symplectic Lie algebroid one has to be careful as to 
what is a formal deformation. One can define the algebra of $E$-differential
operators ${}^E Op (M)$ as the abstract algebra generated by $C^{\infty} (M)$ 
and $\Gamma (M,E)$ subject to obvious relations. This is in itself a 
deformation of the commutative algebra $\Gamma (M,S^* E)$. The morphism $\rho$ 
defines a (not necessarily injective) homomorphism of ${}^E Op (M)$ to the algebra 
of all differential 
operators on $M$. It is therefore more natural 
to call a deformation a formal power series
$$
\varpi =\sum_{k \geq 0} (i\hbar )^k P_k
$$
where $P_k$ are $E$-bidifferential operators, $P_0 (f,g)=fg$ and 
$P_1 (f,g)-P_1 (g,f) = \{ f, g\} $. One imposes a natural associativity condition on $\varpi$.
 Isomorphisms of deformations are defined similarly, as well as derivations.

\medskip

The main classification results of this paper can be summarized as follows.

\begin{th}
Let $(E,[\ ,\ ],\rho ,\omega ) $ be a symplectic Lie algebroid on $M$.
The set of isomorphism classes of $E$-deformations of $(E,[\ ,\ ],\rho ,\omega 
) $ is in
bijective correspondence with the space 
$$
 \frac{1}{i\hbar}\omega + \  {}^E \! H^2 (M,\C [[\hbar 
]]) ,
$$
where $ {}^E \! H^2 (M,\C [[\hbar ]])$ is the second cohomology group of 
the $E$-de Rham complex. The
cohomology class $\theta$ associated to the deformation by the above theorem
is called its characteristic class. 
\end{th}

\begin{th}
Let $\App (M)=(C^{\infty} (M)[[\hbar ]]),*)$ be an $E$-deformation of $M$. There exists a Lie algebra
extension 
$$
0\rightarrow \mbox{Ad}(\App (M) )\rightarrow {}^E Der (\App (M) )\rightarrow
{}^E H^1 (M,\C [[\hbar ]] )\rightarrow 0.
$$
\end{th}

In case of complex manifolds (cf. section  \ref{complex symplectic manifolds} which can be read independently of the rest of the paper), a question related to the classification problem above is to classify all deformations of the sheaf of algebras of holomorphic functions on a complex manifold $M$ with a holomorphic complex structure. By definition, such a deformation is a structre of a sheaf of algebras on ${\cal {O}}_M [[\hbar]]$ which is equal to ${\cal {O}}_M $ modulo $\hbar$ and such that the local multiplication law and the transition isomorphisms are given by the power series in $\hbar$ with coefficients in (bi)differential operators. For any deformation, one can define its characteristic class $\theta$ as in the smooth case. Given a complex
manifold $M$, let $F^* \Omega^{(*,*)}$ denote
the decreasing filtration of the de Rham complex of $M$ given by
\begin{equation}
  \label{eq:inthodge}
F^i \Omega^{(*,*)}(M)=\sum_{k\geq i ,l\geq 0} \Omega^{(k,l)}(M).
\end{equation}

\begin{th}
Let $(M,\omega )$ be a complex manifold with a holomorphic symplectic structure $\omega$, such that the 
maps 
\begin{equation}
  \label{complex21}
 H^i (M ,\C )\rightarrow H^i(M,{\cal O}_M) ;\ i=1,2
\end{equation}
are surjective. The set of isomorphism classes of formal deformations of $(M,\omega
)$ is in bijective correspondence with $H^2 (F^1 \Omega^{*,*}
(M),d)[[\hbar ]]$. Moreover there exists a family of smooth (non-affine) maps:
\begin{equation}
  \label{cc10}
\tau_n : \{ H^2 (F^1 \Omega^{*,*} (M),d) \}^n \rightarrow H^2(M,{\cal O}_M) )
\end{equation}
such that the characteristic class of the deformation associated to $\alpha_0
+\hbar \alpha_1 +\ldots $ is given by the sum
\begin{equation}
  \label{eq:intcomp 2}
  \frac{1}{i\hbar}\omega + \sum_n \hbar^n (\alpha_n + \tau_n (\alpha_0 , \ldots , \alpha_{n-1} )).
\end{equation}
The associated formal deformation of the sheaf of algebras of holomorphic
functions is locally isomorphic to the Weyl deformation of holomorphic
functions on an open polydisc in $\C^{dim M}$.
\end{th}

The maps $\tau _n$ are related to Rozanski-Witten invariants and to the homotopy Lie algebra structure of Kapranov \cite{RW}, \cite{Kapranov}. This relation will be discussed in a subsequent paper.

The basic tool used below is a notion of $E$-jets, a generalization of the
notion of infinite jets of smooth functions to the case when the Lie
algebroid in question is not identical with the sheaf of vector fields on
$M$ (comp. \cite{E1}, \cite{E2}, \cite{Rinehart}). In the case when the structural 
morphism $\rho
$ is injective at the level of sections of $E$, the bundle ${}^E Jets$ is the
bundle of algebras of formal Taylor coefficients in the directions given by
vector fields from $\rho (E)$. 

For example, in the case when 
$$
M=\{ x\in \R^n | x_1 \geq 0 ,\ldots ,x_k \geq 0 \},
$$
 $E$-jets at the points of the stratum 
$$
\{ x | x_i=0 \iff i\in \{ i_1 ,\ldots i_l \} \ \}
$$
are formal power series in
$$
x_i, \ i\not{\in } \{ i_1 ,\ldots,i_l \} \ \mbox{and} \ 
                                     log(x_{i_1}),\ldots ,log(x_{i_l}).
$$

The bundle ${}^E Jets$ carries a {\bf natural flat $E$-connection} which we call 
Grothendieck connection (\cite{Gr},\cite{Kash}), and the algebra of horizontal
sections can be identified with the algebra of smooth functions on
$M$. Moreover, the symplectic structure on the Lie algebroid induces on ${}^E
Jets$ a structure of a bundle of Poisson algebras.
 
A formal deformation of $M$ in our general sense can be also defined as 
a fiberwise formal deformation of this bundle of (local, complete) Poisson 
algebras compatible with the Grothendieck connection.

On the other hand, given a symplectic structure $\omega $ on $E$, one can 
construct the associated Weyl bundle $ {}^E \bf W $ whose fibre at a point is 
linearly isomorphic to the completed symmetric algebra 
$\hat{S}^* (E^* )[[\hbar ]]$. This is naturally a bundle of algebras with the 
fibers endowed with the algebra structure given by the Weyl product. 

Now, given a formal deformation $*$, ${}^E Jets $ becomes a bundle of algebras and
one can show that
$$
{}^E Jets \stackrel{\sim}{\rightarrow} {}^E {\bf W }
$$
(non-canonically) as bundles of algebras. The image of the Grothendieck connection under this 
isomorphism provides a flat connection on ${}^E \bf W$ preserving the 
multiplication and having some additional properties (a {\bf 
Fedosov connection}). Thus any formal deformation gives a Fedosov connection
on the Weyl bundle. On the other hand, given a Fedosov connection on the Weyl 
bundle, the algebra of flat sections is linearly isomorphic to 
$C^{\infty}(M)[[\hbar ]]$ and one can define the deformed product $*$ on the 
bundle of $E$-jets.  

Note that a Fedosov connection $\nabla$ is a flat connection with values in the 
Lie algebra ${\frak g} =Der ( {\mathbb A}^{\hbar } )$, where $  {\mathbb A}^{\hbar 
}$ 
denotes the Weyl algebra of the standard linear symplectic space. There is the natural central extension 
$$ 0 \rightarrow {\mathbb C}[[\hbar]] \rightarrow \tilde{\frak
  g} \rightarrow {\frak
  g} \rightarrow 0
 $$

The characteristic class of the deformation can be computed as
the curvature $\frac{1}{2}[ \tilde{\nabla} ,\tilde{\nabla} ]$ for any $\tilde{\frak
  g}$-lifting $\tilde{\nabla}$ of $\nabla$.  

\medskip
 
The index theorems for elliptic operators in various contexts can be obtained as applications of the classification results of this paper. The main point 
is that symbol calculus of pseudodifferential operators is a formal 
deformation in a disguise (see \cite{NTa/a}). Moreover, the standard trace 
on the trace class pseudodifferential operators corresponds to the unique 
trace on a formal deformation of any symplectic manifold.  
As the result, the precise information about formal deformations of a constant 
rank Poisson structure associated to a foliation allows one to prove 
higher index theorems for foliation algebras (see \cite{NTfamilies}). 

More generally, let $E$ be a symplectic Lie algebroid on a manifold $M$. Given its deformation, we define a trace density map 
$$
\mu^{\hbar}: CC_*^{per}(\fA (M)) \rightarrow ({}^E\Omega ^{2n-*}(M)((\hbar))[u^{-1},u]], d)
$$
where the left hand side is the periodic cyclic complex of the deformed algebra of functions. We compute (theorem \ref{index thm for E}) the action of this map on cohomology in terms of reduction modulo $\hbar$ (the principal symbol map)
$$CC_*^{per}(\fA (M)) \rightarrow CC_*^{per}(C^{\infty} (M)),$$
of Connes' morphism $CC_*^{per}(C^{\infty} (M)) \rightarrow \Omega^*(M)[u^{-1}, u]],$ of the characteristic class $\theta$ of the deformation, and of the $\hat{A}$ class of the bundle $E$. This theorem generalizes the index theorems from \cite{19 Fedosov's book}, \cite{NT1}, \cite{NTfamilies}. It allows to give a new proof of a recent theorem of Epstein and Melrose \cite{Epstein Melrose}.
An analogous theorem for complex manifolds, proven in \cite{BNT}, implies a Riemann-Roch theorem for elliptic pairs conjectured by Schapira and Schneiders \cite{SS}.

Note that, conjecturally, theorem \ref{index thm for E} holds for any Poisson structure. This conjecture relies on a general formality conjecture for chains \cite{foch} which is an analog of the formality theorem of Kontsevich \cite{Konts} for the complexes of cyclic chains. If true, the generalized theorem \ref{index thm for E} allows to generalize the $\hat{A}$ class of the tangent bundle of the foliation of symplectic leaves to the case when a Poisson structure is not necessarily regular.

\begin{remark}
\begin{em}
The whole idea of applying formal methods to geometry stems from the papers of
I.M. Gelfand and his collaborators (\cite{Gelfand},\cite{GGL}), and especially one of the aims of
developing the theory of formal deformations and associated characteristic
classes as in \cite{NTfamilies} is to apply the machinery of index theorems in
the case of corners to combinatorial
study of characteristic classes. The second autor is grateful to B.L.Feigin and to I.M. Gelfand
for many fruitful insights of this subject and of its applications to index theory.
\end{em}
\end{remark}

\section{Generalized jet spaces}

\subsection{ Lie algebroids}
We will recall below some of the standard notions connected with Lie algebroids 
\cite{?}, \cite{pradines}, \cite{bb}.

\medskip

\begin{definition}
\begin{em}
Let $M$ be a smooth manifold. A {\bf Lie algebroid} on $M$ is a triple $(E,
\rho , [\ ,\ ])$, where $E$ is a vector bundle $E$ on $M$, $[\ ,\ ]$ is a Lie
algebra structure  on the sheaf of sections of $E$, and $\rho $ is a map of vector bundles 
$$
\rho :E \rightarrow TM 
$$
such that the induced map
$$
\Gamma (\rho ):\Gamma (M,E) \rightarrow Vect(M)
$$
is a Lie algebra homomorphism and, for any sections $\sigma $ and $\tau$ of
$E$ and a smooth function $f$ on $M$, the following identity holds:
$$
[\sigma ,f \tau ] = \rho (\sigma ) (f) \cdot \tau + f [\sigma ,\tau ].
$$
As a matter of notation, we will use
${\frak g}_E$ to denote the Lie algebra $(\Gamma (M,E),[\ ,\ ])$, and will
regard $C^{\infty }(M)$ as a left ${\frak g}_E$-module with the action given by
$\Gamma (\rho )$ and, for $\sigma \in {\frak g}_E $ and $f \in C^{\infty }(M)$,
$$
\sigma f \stackrel{def}{=} (\Gamma (\rho )(\sigma )) f
$$
Note also that ${\frak g}_E$ is a left $C^{\infty }(M)$-module. From now on we 
will abbreviate $\Gamma (\rho )$ by $\rho$. 
\end{em}
\end{definition}

The following construction gives a natural generalization of the de Rham
complex for a Lie algebroid.

\begin{definition}
Let $(E, \rho , [\ ,\ ])$ be a Lie algebroid on $M$. The $E$-de Rham
complex $({}^{E} \Omega^* (M),{}^E d)$ is given by
\begin{equation}
\begin{array}{l}
{}^{E} \Omega^* (M)=\Gamma (M,\Lambda^{\bullet} ( {E}^*)),\\
{}^{E}{\rm d}\omega (\sigma_1,\ldots ,\sigma_{k+1} ) =\\
=\left\{ \begin{array}{rl} 
           {}&  \sum_{i} {(-1)}^i \sigma_i \omega (\sigma_1,\ldots
,\hat{\sigma}_i ,\ldots  ,\sigma_{k+1}) \\ 
  +& \sum_{i< j} {(-1)}^{i+j-1} \omega ([\sigma_i ,\sigma_j ],\sigma_1,\ldots 
,\hat{\sigma}_i
,\ldots ,,\hat{\sigma}_j ,\ldots ,\sigma_{k+1}). 
\end{array} \right\}.
\end{array}
\end{equation}
The cohomology of this complex will be denoted by ${}^E H^* (M)$ and called the
E-de Rham cohomology of $M$.

A $E$-connection on a vector bundle $F$ on $M$ is a linear map
$$
\nabla : \Gamma ( F \otimes \Lambda^{\bullet} ( {E}^*)) \rightarrow \Gamma ( F
\otimes \Lambda^{\cdot +1} ( {E}^*)) 
$$
satisfying the Leibnitz identity:
$$
\nabla (f \sigma )= {}^E df \sigma + f \nabla (\sigma ).
$$
\end{definition}

\medskip

Note that the above definition makes sense because local sections of $E$ are
closed under Lie bracket. However, in
distinction to the standard de Rham complex, $({}^{E} \Omega^* (M),
{}^{E}{\rm d} )$ is not locally acyclic and hence does not give an acyclic
resolution of $C^{\infty } (M)$.

\subsection{E-Differential operators and E-jets}


\begin{definition}

Let $(E, \rho , [\ ,\ ] )$ be a Lie algebroid on $M$. The sheaf of {\bf
E-differential operators} ${ {}^E Op} $ on $M$ is given by
$$
U \mapsto { {}^E Op}(U) = T /   \begin{array}{c}
                                                          \\ 
                                                     \left\{
                                       \begin{array}{c}
                           \sigma \tau - \tau \sigma -[\sigma ,\tau] \\
                          \sigma (f\tau ) -(f \sigma )\tau -\sigma (f)\tau \\
                           f\cdot \sigma - f \sigma
                                       \end{array}
                                                          \right\}.

                                           \end{array}
$$
Here $T$ is the graded algebra freely generated by the algebras  $C^{\infty } (M)$ 
(of
degree 0) and  the free algebra of ${{\frak g}_E} $ (of degree 1). The sheaf of algebras ${ {}^E Op}
$ is equipped with the left action of the sheaf $\cal O_M$ of smooth functions on 
$M$.
Moreover, the grading of the algebra $ T$ descends to a
filtration ${ {}^E Op}_n$ of ${ {}^E Op}$. The sections of
${ {}^E Op}_n$ will be called E-differential operators on $M$ of degree $\leq n$.

The subsheaf of $E$-differential operators without zero order term,
${ {}^E Op}^+ (M) $ is defined as the left ideal generated by the image of the map

$$
E \rightarrow{ {}^E Op}(M).
$$
\end{definition}

\medskip

Let us record the following result 
\begin{lemma}[Poincar\'{e}-Birkhoff-Witt]
$$
Gr { {}^E Op}(M) \simeq \Gamma (M,S(E))
$$
\end{lemma}
We will refer for the proof to \cite{MN}.

\medskip

\begin{remark}
\begin{em}
Note that in the degenerate case, when $\rho$ vanishes, this is the standard
Poincar\'{e}-Birkhoff-Witt theorem applied fiberwise to a bundle of Lie algebras.
\end{em}
\end{remark}

\begin{definition}
Let $(E, \rho , [\ ,\ ])$ be a Lie algebroid on $M$. The space of {\bf E-jets}
on $M$ is the linear space
$$
{}^E \underline{Jets}(M) 
                      = Hom_{C^{\infty }(M)} ({{}^E Op}(M) ,C^{\infty }(M)).
$$
We set
$$
{}^E \underline{Jets}^+(M) = 
                     Hom_{C^{\infty }(M)} ({{}^E Op}^+ (M) ,C^{\infty }(M)). 
$$
\end{definition}
\medskip

\begin{proposition} 
In the above notation, ${}^E \underline{Jets}(M)$ is the space of global
sections of a profinite-dimensional vector bundle ${}^E Jets$. The equation
$$
 \forall \ D\in {{}^E Op (M)},\ \sigma \in {\frak g}_E, \
               \nabla_G ({\sigma } )l (D )  = \sigma  l (D )-l(\sigma D)\ .
$$
defines on the bundle ${}^E Jets$ a flat $E$-connection which we will call the 
{\bf Grothendieck connection}.
\end{proposition}
{\em Proof}. Let $U$ be a local coordinate system on $M$ such that $E|U \simeq
U\times \R^n $. We denote by ($e_1 ,\ldots ,e_n$) the associated basis for
$\Gamma (U;E)$. By Poincar\'{e}-Birkhoff-Witt theorem,
\begin{equation}
  \label{E-operators}
e_{\alpha } =\prod_{i=1,\ldots ,n} \frac{  e_i^{\alpha_i}
  }{\alpha_i !}, \alpha =(\alpha_1 ,\ldots ,\alpha_n )\in ({\mathbb N}\cup
0)^n 
\end{equation}
form a basis of E-differential operators over $C^{\infty } (U)$. Let
$l_{\alpha}$ denote the family of E-jets given by 
$$
\ l_{\alpha } (e_{\beta })=\delta_{\alpha ,\beta}.
$$
The $C^{\infty } (U)$-linear map
$$
C^{\infty} (U)\otimes  \hat{S}\R^n  \simeq  {}^E \underline{Jets} (U).
$$
given by sending the symmetric tensor 
$$
 \delta_1^{ \alpha_1} \ldots 
\delta_n^{ \alpha_n}  ; \ \{ \delta_i \}_{i=1,\ldots ,n}\ \mbox{the
  standard basis of $\R^n$}
$$
to $l_{\alpha}$ defines local trivializations and hence a bundle structure on
the E-jets with the fiber $\hat{S}\R^n$. 

The transition functions of this bundle are given by symmetric powers of matrices with
smooth coefficients, multiplied my smooth functions from an intersection of two open sets to the pronilpotent group $G^{\geq 1}$ of continuous automorphisms of $\hat{S}\R^n$ whose linear part is the identity. Thus the above construction does give a profinite-dimensional vector
bundle. The fact that $\nabla_G$ is a flat $E$-connection is a direct
reformulation of the definition of a Lie algebroid.

\medskip

\begin{example}
\begin{em}
Suppose that $M$ is just a single point. Then the Lie algebroid is given by a
Lie algebra $\g$, the E-jet bundle coincides with the completed symmetric
algebra $\hat{S}\g^*$ and the E-de Rham complex with coefficients in
${}^E$Jets with the Grothendieck connection gives  
$$
(\Lambda^{\bullet} \g^* \otimes \hat{S}^{\bullet} \g^* ,\partial )
$$
where $\partial$ denotes the Koszul differential.
\end{em} 
\end{example}

\medskip

\begin{proposition}
The map
$$
\begin{array}{ccc}
C^{\infty} (M) & \stackrel{{}^E J}{\rightarrow }  & {}^E  Jets (M) \\
       f       & \mapsto                       & \{ D \mapsto \rho (D)f \} 
\end{array}
$$
is an isomorphism onto the kernel of the Grothendieck connection $\nabla_G$.
\end{proposition}
{\it Proof}. Given an element $l$ of $ {}^E \underline{ Jets} (M)$, $\nabla_G  l 
=0$ 
is, by flatness of  $\nabla_G $, equivalent to
$$
l(D) = \rho (D) l(1).
$$
This implies the claim of the proposition.

\medskip

Let $\Delta_0$ be the coproduct on the universal enveloping algebra $U({\frak
g}_E)$ of the Lie algebra ${\frak g}_E$. It is a homomorphism of algebras
\begin{equation}
\begin{array}{cccc}
                      & U({\frak g}_E)& \rightarrow &U({\frak g}_E) \otimes 
                                                          U({\frak g}_E), \\
\forall{\sigma \in {\frak g}_E }: & \sigma & \stackrel{\Delta_0}{\mapsto} & 
                                  \sigma \otimes 1 + 1 \otimes  \sigma . 
\end{array}
\end{equation}

\begin{proposition}
The dual of the coproduct $\Delta_0$ induces on both ${}^E Jets$
and ${}^E Jets^+ $ the structure of bundles of commutative algebras with
fibers at a point $m$ given by  $\hat{S} (E_m^* )$  (respectively the
augmentation ideal $\hat{S}^+
(E_m^* ) $). The Grothendieck connection is a derivation with respect
to this algebra structure. The bundles of algebras ${}^E Jets$ and $\hat{S}(E^*)$ are (non-canonically) isomorphic.

\end{proposition}
{\it Proof}.
To begin with, note that an $E$-jet is a linear map
$$
l\in Hom(U({\frak g}_E),C^{\infty }(M)) 
$$
such that
$$
\begin{array}{c}
\forall {\sigma \in{{\frak g}_E}}:\  l(f\sigma \ldots )=fl(\sigma \ldots ), \\
\forall {\sigma ,\tau \in {\frak g}_E}:\  l(\ldots \sigma (f\tau )\ldots
)=l(\ldots (f\sigma ) \tau \ldots )+l(\ldots \rho (\sigma )(f) \tau \ldots ) 
\end{array}
$$
From this it is easy to check that the transpose of $\Delta_0$, given by
$$
(l_1 l_2 )(D) =( l_1 \otimes l_2 )(\Delta_0 (D)),
$$ 
is in fact well defined and, since $\Delta_0$ is symmetric, defines a 
commutative algebra structure on $E$-jets. 
Since all the other statements are local, it suffices to work within a
trivializing neighbourhood $U$ of a point of $M$ and we will use the
representation (\ref{E-operators}) for the elements of ${{}^E Op}(U)$.
The filtration by order on ${{}^E Op}(U)$ induces a complete
decreasing filtration 
$$
{}^E \underline{Jets} (U) = F_0 \supset F_1 \supset F_2 \supset \ldots
$$ 
on ${}^E \underline{Jets} (U)$ and the Poincar\'{e}-Birkhoff-Witt theorem implies
that 
$$
Gr ({}^E \underline{Jets} (U))=C^{\infty} (U) \otimes \hat{S}({\R^n}^* ) . 
$$
Since the complete local algebra $\hat{S}({\R^n}^* )$ has no deformations in the
class of commutative algebras, this implies in particular that
$$
{}^E \underline{Jets} (U) \simeq C^{\infty} (U) \otimes \hat{S}({\R^n}^* ) .
$$
An explicit isomorphism can be constructed as follows. Let $ (l_1 , \ldots
,l_n )$ be the $E$-jets defined by 
$$
l_i (e_j )=\delta_{ij} , \ l_i (e_{k_1} e_{k_2} \ldots e_{k_l} )=0 \ \mbox{for
                                                                      $l>0$}. 
$$
It is easy to see that the map
$$
{}^E \underline{Jets} (U) \ni l \mapsto \sum_{\alpha } l(e_{\alpha}) e^{\alpha}
$$
gives the required isomorphism ($\{e^{\alpha}\}$ being the dual basis to $\{e_{\alpha}\}$).

To prove the last statement of the proposition, note that for a good cover $\{U_i\}$ one has the isomorphisms $\phi_i: {}^E Jets | U_i \simeq \hat{S}(E^*) | U_i$. The transition isomorphisms $g_{ij} = \phi _i \phi _j ^{-1}$ take values in the pronilpotent group of those automorphisms of $\hat{S}(E*)$ whose derivative at zero is equal to the identity. Therefore the cocycle $\{ g_{ij}\}$ is cohomologous to the trivial cocycle.


\subsection{ Symplectic Lie algebroids}

\begin{definition}
\begin{em}

A {\bf symplectic Lie algebroid} structure on $M$ is a pair $((E, \rho , [\ ,\
]),\omega )$, where $(E, \rho , [\ ,\ ])$  is a Lie algebroid structure on $M$
and $\omega $ is a closed ${E} $-two-form on $M$ such that the associated
linear map: 
$$
{\frak g}_E \times {\frak g}_E \ni (X,Y) \mapsto \omega (X,Y) \in C^{\infty}
(M) 
$$
is a symplectic structure on $E$.
Whenever possible, we will abreviate the notation for the symplectic Lie
algebroid to $(E,\omega )$. 
\end{em}
\end{definition}

\bigskip

A symplectic Lie algebroid structure $(E,\omega )$ on $M$ leads to an
associated Poisson bracket which we will construct and study below. To begin
with, the definition above gives us an isomorphism: 
$$
I_{\omega }: E_m \rightarrow E_m^* .
$$

Given a smooth function $f$ on $M$, we define the associated {\bf Hamiltonian
vector field} $H_f$ as the image of $f$ under the composition:
\begin{equation}
\begin{array}{c}
H: C^{\infty} (M) 
\stackrel{\rm d }{\rightarrow} \Omega^1 (M)
\stackrel{\rho^t }{\rightarrow} \Gamma (M,E^* )
\stackrel{I_{\omega }^{-1} }{\rightarrow} \Gamma (M,E )
\stackrel{\rho }{\rightarrow} Vect(M).     
\end{array}
\end{equation}

\begin{lemma}
\begin{em}
Let $(E,\omega )$ be a symplectic Lie algebroid structure on $M$. The
equation
\begin{equation}
\{ f ,g \} =H_f g
\end{equation}
defines a Poisson structure on $M$. Moreover the following identity holds
\begin{equation}
[H_f ,H_g ]=H_{ \{ f,g \} }.
\end{equation}
  
\end{em}
\end{lemma}
Proof.
It is obvious from the construction that $H_f$ is a vector field in
$\rho ({\frak g}_E )$. Also from the construction we get the equality:
$$
\{ f, g \} = \omega (I_{\omega }^{-1}\rho^t df ,I_{\omega
}^{-1}\rho^t dg ) .
$$
This gives the skew symmetry of $\{ f, g \}$, while the (Jacobi) identity:
$$
\{ f, \{ g, h\} \} + \{ g, \{ h, f\} \} +\{ h, \{ f, g\} \} =0 
$$
is equivalent to 
$$
{}^E d(I_{\omega }^{-1}\rho^t df,I_{\omega }^{-1}\rho^t dg,I_{\omega }^{-1}\rho^t 
dh ) =0
$$
(recall that $\omega$ is a closed E-form).

\medskip

\begin{definition}[Poisson bracket on jets]
\begin{em}
Let, as above, $(E,\omega )$ be a symplectic Lie algebroid structure on
$M$. Let 
$$
\varpi = \sum_i e_i \otimes f_i 
$$
be the antisymmetric tensor in $ \Gamma (M,E^{\otimes 2})$ which is the image
of $\omega $ under the isomorphism $I_{\omega} \otimes I_{\omega}$. The $\{ \ ,\ 
\}$ is the skew-symmetric $C^{\infty}
(M)$-bilinear map given by 
\begin{equation}
\begin{array}{ccc}
 {}^E \underline{Jets} (M)\times {}^E \underline{Jets} (M) &  \stackrel{\{ \
,\ \}}{\longrightarrow } & 
                                                           {}^E
\underline{Jets} (M)   \\ 
      ( l_1 , l_2) &  \mapsto  & \{ D\mapsto (l_1 \otimes l_2 ) ( \Delta_0 (D)
\cdot \sum_i e_i \otimes f_i )\} .
\end{array}
\end{equation}
In the case of $\rho =0$ this is known as the bracket of Berezin-Kirillov-Kostant-Souriau.   
\end{em}
\end{definition}

\medskip

\begin{th}
(${}^E {Jets} ,\{ \ ,\ \} $) is a bundle of Poisson algebras over
$M$, with fiber isomorphic to the Poisson algebra
$$
(\R [[ x_1 ,\ldots ,x_n ,\xi_1 , \ldots , \xi_n ]] , \{ \ ,\ \}_{st} )
$$
where n=$\frac{1}{2}$rank$E$ and
$$
\{ f ,g \}_{st} = \sum_{i=1}^{n} \left( \partial_{x_i} f\partial_{\xi_i} g
                     -\partial_{\xi_i} f\partial_{x_i} g \right)
$$ 
\end{th}
Proof.
To begin with, let us prove that (${}^E {Jets} ,\{ \ ,\ \}
$) is a bundle of Poisson algebras. Since this is a local statement and since
$\omega $ is nondegenerate, we can assume that ($e_1, \ldots ,e_n ,f_1, \ldots
,f_n$) form a basis for the space of sections of $E$. We set 
$$
\varpi = \sum_i e_i \otimes f_i .
$$
From the definition it follows immediately that the following equivalences
hold:
$$
\begin{array}{lccc}
0)&
  \left\{
          \begin{array}{c}
  \forall_{l_1 ,l_2 ,l_3} \ \{ l_1 ,l_2 \} =- \{ l_2 ,l_1 \} \\
      \end{array} 
  \right\}  
           &  \Longleftrightarrow &  
                      \left\{  \ \mbox{$\varpi $ is antisymetric} \ 
          
                                                      \right\}   \\
   &&&\\
1)&
  \left\{
          \begin{array}{c}
  \forall_{l_1 ,l_2 ,l_3} \ \{ l_1 ,l_2 l_3 \} =\\
      l_2 \{ l_1 ,l_3 \} +l_3 \{ l_1 ,l_2 \}
          \end{array} 
  \right\}  &  \Longleftrightarrow &  
                      \left\{
          \begin{array}{c}
                 \iota \otimes \Delta_0 ( \Delta_0  (D) \cdot \varpi )= \\
(id + \sigma_{23} ) ( (\Delta_0 \otimes \iota (\Delta_0
(D)) \cdot \varpi \otimes 1 )
          \end{array}
                                                      \right\}   \\
  &&&\\       
2)&
  \left\{
          \begin{array}{c}
  \forall_{l_1 ,l_2 ,l_3} \  \{  l_1 ,\{ l_2 ,l_3 \} \} \\
     +\mbox{ \ cyclic permutations} =0
          \end{array} 
  \right\}  &  \Longleftrightarrow &  
                      \left\{
          \begin{array}{c} 
  \sum_{ij} e_i e_j \otimes f_i \otimes f_j +e_i \otimes f_i e_j \otimes f_j
\\
  + \mbox{\ cyclic permutations} =0
          \end{array}
                                                      \right\}   \\

\end{array}
$$
But 0) holds by construction, 1) is a straightforward consequence of the
coassociativity of $\Delta_0$ and the way it acts on ${\frak g}_E$, while 2) is
easily seen to be equivalent to  
$$
[\varpi , \varpi ] =0 \ \mbox{\ in \ }\  (\Gamma ( M, \Lambda^3 E ), 
\mbox{Schouten bracket} ),
$$
which is in turn equivalent to the fact that $\omega $ is a closed $E$-form.

To prove the rest of the theorem, it is usefull to get a more explicit
representation of the Poisson structure. So still working locally, let us
recall that
$$
{}^E \underline{Jets} (U) \simeq C^{\infty } (U) \otimes \hat{S} ({\R^{2n}}^*),
$$
with $l \in \Gamma(U;E^* )$ giving a generating set of first order $E$-jets
(in the grading of $\hat{S}$). Since $\omega$ is non degenerate, we can choose
as first order sections 
$$
e^* = \omega (e, \cdot ), \ \mbox{\ for} \ e\in \Gamma(U;E ).
$$ 
Now
$$
\{ e^* , f^* \} (1) = \omega (e,f) 
$$
by the definition of the bracket, and hence
$$
\{ e^* , f^* \} = \omega (e,f) 1 + \ \mbox{\ higher order terms } .
$$
In other words, $\{ \ , \ \}$ gives a Poisson structure on $\R [[ x_1 ,\ldots
,x_n ,\xi_1 , \ldots , \xi_n ]]$ associated to a symplectic structure. An
application of the formal Darboux theorem finishes the proof.

\medskip


\section{Formal deformations}


\subsection{Basic definitions}
\begin{definition}\label{5.1.1}
\begin{em}
Let $A$ be an associative unital algebra over a unital ring $k$. A formal deformation of $A$
is a structure of an associative algebra over $k[[\hbar ]]$ on $A[[\hbar ]]$
given by a product $*$ of the form
\begin{equation}
f*g =fg + \sum_{i=1}^{\infty} {\hbar}^i \varpi_i (f,g)
\end{equation}
and
$$ 1*f=f*1=f$$
An isomorphism of two deformations $*$ and $*'$ is a formal series $T(a) = 1 + \sum_{k \geq 1} \hbar ^k T_k (a)$ such that $T(a)*T(b) = T(a*'b)$. A deformation quantization of a smooth manifold $M$ is a deformatin of $C^{\infty}(M)$ for which $\varpi _k$ are bidifferential operators. An isomorphism of such deformation quantizations is an isomorphism of corresponding deformations for which $T_k$ are differential operators.
\end{em}
\end{definition}

\begin{notation}
Given a formal deformation of a Poisson manifold $(M, \{ ,\ \} )$, the algebra
$(C^{\infty} (M)[[\hbar ]],*)$ will be denoted by
$\App (M)$.
\end{notation}

\bigskip
A general construction of deformations of Poisson structures is given in
\cite{Konts} 
\subsection{Weyl deformation}

\begin{definition}
\begin{em}
Let $(V,\omega )$ be a symplectic vector space over a field $k$ containing the
square root $i$ of $-1$. Let $TV$
denote the tensor algebra of $V$. The Weyl algebra of $(V,\omega )$ is the
associative algebra over the ring of formal power series $k[[\hbar ]]$ given by
\begin{equation}
{\W }(V)=TV/(v\otimes w -w\otimes v -i\hbar \omega (v,w)),
\end{equation}
completed in $(\hbar ,V)$-adic topology.
\end{em}
\end{definition}

Note that
$$
V \leadsto \W (V)
$$
is a functor from the category of finite dimensional symplectic vector spaces
to the category of finitely generated complete graded
algebras over $k[[\hbar ]]$. The grading on $\W (V)$ is 
\begin{equation}
|\hbar | =2 \mbox{\ and, for any $v\in V$, \ } |v|=1.
\end{equation}

A particular case of this definition will deserve a separate name. Let
$$
V=\R^{2n}
$$
with coordinates
$$
x=(x_1 ,\ldots ,x_n ), \ \xi =(\xi_1 ,\ldots ,\xi_n ). 
$$
Let $\omega_{st}$ be the symplectic form on $V$ given by
$$
\omega_{st}(x_i ,x_j) = \omega_{st}(\xi_i ,\xi_j )=0,\  \omega_{st}(\xi_i ,x_j
)=\delta_{i,j}. 
$$
The algebra $\App (\C^{2n} ,\omega_{st} \otimes_{\R} 1 )$ will be called the
{\bf Weyl algebra} and will be denoted by $\hat{\mathbb A}^{\hbar}$. Its
generators, image of the above basis for $\R^{2n}$ in $\hat{\mathbb A}^{\hbar}$,
will be denoted by
\begin{equation}
\fx_1 , \fxi_1 , \ldots ,\fx_n , \fxi_n .
\end{equation}

{\bf {Notation.}} 
Denote by $\tilde{\frak g } $ the Lie algebra 

$$
\{ \hbar^{-1} \Phi |\Phi \in \hat{\mathbb A}^{\hbar} , i\Phi \mbox{\ real mod } \
\hbar \} , 
$$

and by $\frak g $ the quotient:
$$
\{ \hbar^{-1} \Phi |\Phi \in \hat{\mathbb A}^{\hbar} , i\Phi \mbox {\ real mod }\
\hbar \} / \{ i\hbar^{-1} \R  +\C [[\hbar ]] \} , 
$$
both with the bracket given by

$$
[f,g]= f*g-g*f.
$$

We give $\tilde{\frak g }$ the grading
$$
\tilde{\frak g } =\prod_{n\geq {-2}} \  {\tilde{\frak g } }_n,
$$
with the grading induced by the grading of $\hat{\mathbb A}^{\hbar}$ and set
$$
{\tilde{\frak g }}^0=\prod_{n\geq 0} \  {\tilde{\frak g } }_n .
$$
We will use the same notation for the induced grading on $\frak g $, so that

$$
{\frak g } =\prod_{n\geq {-1}} \  {{\frak g } }_n
$$

\medskip

Note that the group $G^0 $ of continuous automorphisms of $\hat{\mathbb
A}^{\hbar} $ is a profinite dimensional Lie group with the Lie algebra ${\frak
g}^0  $.
It contains as a subgroup the group $G_0=\mbox{Sp}(2n)$ of linear automorphisms of
$\hat{\mathbb A}^{\hbar} $ and the quotient
$$
G^0/\mbox{ Sp }(2n)
$$
is contractible ($\cong \R^{\infty }$). We set
$$
G^n = {\rm exp }({\frak g }_{\geq n }).
$$

\begin{definition}
\begin{em}
The {\bf Weyl deformation} of $\R^{2n}$ is the formal deformation of the
Poisson manifold $(\R^{2n} ,\omega_{st})$ given by the ({\bf Moyal}) product  
\begin{equation}
(f*g) (x,\xi )={\rm exp}\left( \frac{i\hbar}{2} \sum_{k=1}^{n} \left(
\partial_{x_k} \partial_{{\eta}_k} - \partial_{\xi_k} \partial_{y_k}
\right) \right) f(x,\xi )g(y ,{\eta})\mid_{ (x=y ,\xi
=\eta)} . 
\end{equation}
We denote by $\App_c (\R^{2n})$ the ideal of $(\App (\R^{2n}) ,*)$ consisting
of formal power series in $\hbar $ with coefficients of compact support. $(\App
(\R^{2n}) ,*)$ will be always considered as a topological algebra, with the
$\hbar $-adic topology (and $C^{\infty } $-topology in coefficients).

\end{em}
\end{definition}

\medskip

Since the Moyal product is local, it defines a sheaf of associative algebras
on $\R^{2n} $:
$$
U \rightarrow (\App (U),*),
$$
where $\App (U)$ is the linear space $C^{\infty }(U)[[\hbar ]] $.
 
\begin{th}\label{local trace}
Let $U$ be an open connected subset of $\R^{2n}$. The center of $\App (U)$ is
$\k \cdot 1$ and $\App_c (U)$
has a unique (up to a scalar multiple), continuous, $\C [\hbar^{-1} ,\hbar
]]$-valued trace $Tr$ given by 

\begin{equation}
Trf=\frac{1}{(i\hbar )^n n!} \int_U f\omega^n_{st}
\end{equation}
\end{th}
Proof: cf. \cite{19 Fedosov's book} 
\bigskip

The Weyl deformation is a special case of the following construction. Let the
coordinates on $\R^{2n}$ be given by
$$
x= (z_1 ,\ldots ,z_{2k} ,y_1 ,\ldots , y_l ,\eta_1 ,\ldots ,\eta_l , x_1,\ldots
,x_m ,\xi_1 ,\ldots ,\xi_m ) 
$$
and set
\begin{equation}\label{standard bracket}
\{ \ ,\ \} =\sum_{i\leq k} z_i z_{i+k}\partial_{z_i}\wedge \partial_{z_{i+k}} +
\sum_{i\leq l} y_i \partial_{y_i}\wedge \partial_{\eta_{i}} +
\sum_{i\leq m} \partial_{x_i}\wedge \partial_{\xi_{i}}.
\end{equation}

\medskip

\begin{definition}
\begin{em}
The Weyl deformation of $\R^{2n}$ associated to the Poisson structure
(\ref{standard bracket}) is given by 
$$
f*g = {\rm exp}\left( \frac{i\hbar}{2} \sum_{i=1}^{n} 
D_i \otimes E_i -E_i \otimes D_i  \right) f(x )g(\bar{x} )\mid_{ x=\bar{x}} , 
$$
where
$$
D_i = \left\{ \begin{array}{ll}
               z_i \partial_{z_i} & \mbox{for $i=1, \ldots , k$} \\
               y_i \partial_{y_i} & \mbox{for $i=2k+1, \ldots , 2k+l$} \\
               \partial_{x_i} & \mbox{for $i=2k+2l+1, \ldots , 2k+2l+m$} 
              \end{array}  
                                         \right. ,
$$
and
$$
E_i = \left\{ \begin{array}{ll}
                z_{i+k} \partial_{z_{i+k}} & \mbox{for $i=1, \ldots , k$} \\
               \partial_{\eta_i} & \mbox{for $i=k+1, \ldots , k+l$} \\
               \partial_{\xi_i} & \mbox{for $i=k+l+1,\ldots ,k+l+m$} 
               \end{array}  
                                           \right.
$$
\end{em}
\end{definition}

\subsection{Formal deformations associated to symplectic Lie algebroids, Fedosov
construction }

Let $(E,[,],\rho ,\omega )$ be a symplectic Lie algebroid over a smooth
manifold $M$. Recall that we have associated to $E$ the following structures.
\begin{itemize}
\item A Poisson structure $\{ \
,\ \}_E $ on $M$ given by a skew-symmetric tensor 
$$
\varpi_1 \in {\frak g}_E \otimes_{C^{\infty }(M)} {\frak g}_E.
$$
\item A left ${\cal O}_M$-module ${}^E {Op}$ (the sheaf of E-differential
  operators).
\item The bundle ${}^E Jets$ of Poisson algebras isomorphic (not canonically) as a profinite
vector bundle to $\hat{S}(E^*)$.
\end{itemize}

 We set $n=\frac{1}{2}$rank$ E$ and will {\bf fix} this
notation throughout this section. 

\medskip

\begin{definition}
\begin{em}
An $E$-deformation of $M$ is a formal deformation of the Poisson manifold
$(M,\{ \ ,\ \}_E )$ with a $*$-product of the form
\begin{equation}
f*g = f\cdot g +\sum_{l,k\geq 1} \left((i\hbar
  )^k 
\rho(D_{l,k} )(f) \cdot \rho ( E_{l,k} )(g)\right) ,  
\end{equation}
where $D_k$ and $E_k$ are E-differential operators on $M$, the tensor
$$
\varpi_E =1\otimes 1 + \frac{1}{2}i\hbar \varpi_1 +\sum_{l,k\geq 2} \left((i\hbar )^k
(D_{l,k} ) \otimes ( E_{l,k} )\right) \in  {}^E {Op}(M)
\otimes_{C^{\infty }(M)} {}^E {Op}(M) 
$$
satisfies the equation
\begin{equation}\label{assoc}
(\Delta_0 \otimes id) (\varpi_E ) \cdot 1\otimes \varpi_E =(id\otimes
\Delta_0) (\varpi_E ) \cdot \varpi_E \otimes 1 ,
\end{equation}
and 
$$ \frac {1}{i \hbar}(f*g - g*f ) = \{f, g\}_E + O(\hbar) $$
The corresponding associative algebra $( C^{\infty} (M)[[\hbar ]] ,*)$ will be
denoted by $\App (M)$.
\end{em}
\end{definition}

\medskip

This definition has as a corollary the following lemma, which will allow us to
think of the {\em category} of $E$-deformations. 

\begin{lemma}
Let $(E,[,],\rho ,\omega )$ be a symplectic Lie algebroid
on $M$. Given an $E$-deformation, the associated tensor $\varpi_E$ induces a
structure of a profinite dimensional bundle of associative algebras on
$$
{}^E Jets \otimes_{\R} \C [[\hbar ]] .
$$
We will denote this bundle of algebras by $({}^E Jets  ,* )$. 
The Grothendieck connection extends to a flat connection $\nabla_G \otimes_\C
id $ satisfying the Leibnitz identity: 
$$
\nabla_G (l_1 *l_2 ) =\nabla_G (l_1 ) * l_2 + l_1 * \nabla_G (l_2 ).
$$

\end{lemma}
{\em{Proof}}. The definition of a $E$-deformation gives a tensor of the
form  
\begin{equation}
\varpi_E =\sum_{i_k ,k} \left((i\hbar )^k D_k^{i_k}\otimes_{C^{\infty }(M)} 
 E_k^{i_k} )\right) ,
\end{equation}
where $D_k^{i_k}$ and $D_k^{i_k}$ are $E$-differential operators on $M$. We
set, for a pair of sections $(l_1 ,l_2)$ of ${}^E Jets \otimes_{\R} \C [[\hbar
]]$,
\begin{equation}
(l_1 *l_2) (D)=(l_1 \otimes l_2) (\Delta_0 (D) \cdot \varpi_E ).
\end{equation}
Let $X \in \Gamma (M, E)$ and $E,F$ be $E$-differential operators on
$M$. The identity
$$
\Delta_0 (X) (fE\otimes F - E \otimes fF)=(Xf \otimes 1 - 1 \otimes Xf)
E\otimes F
$$ 
shows that $l_1 *l_2$ depends only on the class of $\varpi_E $ in 
$$
{}^E {Op} (M) \otimes_{C^{\infty }(M)} {}^E {Op} (M).
$$
The associativity of the $*$-product on ${}^E Jets$ is equivalent to (\ref{assoc}).
\bigskip

\begin{definition}
\begin{em}
A morphism of two $E$-deformations $*_1$ and $*_2$ of $M$ is an algebra
homomorphism of spaces of sections of the associated $E$-jet bundles
$$
\Phi :( \Gamma (M, {}^E Jets ),*_1 )\rightarrow ( \Gamma (M, {}^E Jets ),*_2 )
$$
which preserves the subspace of $\nabla_G$-flat sections in $\Gamma (M, {}^E
Jets )$. A derivation of an $E$-deformation $*$ is a derivation of the algebra
$(\Gamma (M, {}^E Jets ),*  )$ which preserves the subalgebra of
$\nabla_G$-flat sections.  
\end{em}
\end{definition}

\medskip

\begin{remark}
\begin{em}
The main point of the above choice of definitions is the fact that in our
general context the algebra $(\App (M) ,* )$ does not carry enough information
about the Lie algebroid $E$ to determine the product on $E$-jets. As a typical
example, in the case when the structure map
$$
\rho : E \rightarrow T(M) 
$$
is zero, an $E$-deformation is a nontrivial deformation of a bundle of
symplectic Lie algebras preserving the bundle structure, while $(\App (M) ,*)$
is just the undeformed algebra of $\C [[ \hbar ]]$-valued smooth functions on
$M$. However, in cases of most interest for us $\rho $ will be injective on
${\frak g}_E$ and in this case the deformation of the algebra of smooth
functions for which the *-product is given by $E$-bidifferential operators has a 
unique extension to a deformation of the bundle ${}^E Jets$ and hence defines an $E$-deformation. As it turns out, the replacement of the
algebra of smooth functions by the space of jets makes most of the theory more
transparent.   
\end{em}
\end{remark}

\subsection{$E$-differential forms with coefficients}
Let $\mathbb L$ be a profinite dimensional U(n)-module. Define the space of {\bf $\mathbb L$-valued $E$-differential forms}
\begin{equation}
{}^E \Omega (M, {\mathbb L})
\end{equation}
as follows.
An element of ${}^E \Omega (M, {\mathbb L})$ is a collection $s_U$ of elements
of ${}^E \Omega (U, \Lambda E^*) \otimes {\mathbb L}$ subject to 
$$
s_U = g_{UV} s_V, \ g_{UV} : U\cap V \rightarrow \mbox{U(n)}
$$
where $g_{UV}$ are the transition functions of the bundle $E$ (we reduce the structure group of $E$ to the maximal compact subgroup $U(n)$).

\medskip

\begin{definition}
\begin{em}
Let $\varpi $ be the symplectic form on $E^*$ given by
   \begin{equation}
\varpi (I_{\omega }(v) ,I_{\omega }(w)) =\omega (v,w). 
   \end{equation}
Then $A_{-1}$ is the element of ${}^E \Omega^1 (M, \tilde{\frak g})$ given by
   \begin{equation}
A_{-1 }: E_m \stackrel{\hbar^{-1 }I_{\omega_m }}{\rightarrow} E_m^*
\hookrightarrow \App (E_m^* ,\varpi_m ).
   \end{equation}

\end{em}
\end{definition}

\medskip

\begin{lemma}\label{acyclicity}
The operator $A_{-1}$ satisfies the identity 
   \begin{equation}
[A_{-1},A_{-1}]=(i\hbar)^{-1} \omega 
   \end{equation}
as elements of $ {}^E \Omega
(M,centre (\tilde{\frak g}))\subset {}^E \Omega (M, {\mathbb A}^{\hbar }) $.
In particular $[A_{-1},A_{-1}]$ vanishes in $ {}^E \Omega^2 (M,{\frak
g})$. The associated complex   
$$
( {}^E \Omega (M, \vfA) ,\mbox{ad} A_{-1} )
$$
is acyclic in positive dimension and its zeroth cohomology group coincides with
$ C^{\infty } (M)[[\hbar ]]$.
\end{lemma}
{\em Proof}. The first identity is straightforward. The rest of the statement follows from the fact that $ad \; A_{-1}$ can be identified with the Koszul differential on $\Gamma (M, \wedge E^* \otimes \hat{S}(E^*)[[\hbar]])$. 

\subsection{Fedosov construction}

Let ${\cal {P}}'$ be the bundle of symplectic frames in $E$. Let  ${\cal {P}}$ be a reduction of this principal $Sp(2n)$-bundle to the maximal compact subgroup $U(n).$

\begin{definition}[Weyl bundle, Fedosov connection]
\begin{em}
The bundle
$$ 
 {}^E  {\bf {W}} = {\cal P} \times_{\mbox{U(n)}} \vfA
$$
is called the {\bf Weyl bundle} of $E$. A linear map
$$
\nabla: {}^E \Omega^0 (M, \hat{\mathbb A}^{\hbar })\rightarrow  {}^E \Omega^1 (M,
\hat{\mathbb A}^{\hbar }) 
$$ 
is called a flat connection on the Weyl bundle if it satisfies the equations
$$
\begin{array}{l}
\nabla (vw)=v\nabla (w) + \nabla (v) w \\
\nabla ^2 =0.
\end{array}
$$
It is called a {\bf Fedosov connection} if it is flat and if there exists a
${\frak g}_0$-connection $\nabla_0$ in $E$ such that
$$
\nabla = \nabla_0 + A_{-1} + \sum_{i\geq 1} A_i , \ A_i \in \Omega^1 (M, \g_i),
$$
(recall that ${\frak g}_0 = {\frak sp}$(2n)).
\end{em}
\end{definition}

\bigskip

\begin{th}[Fedosov construction]\label{Fedosov construction}
Let $\theta $ be an element of 
$$
(i\hbar)^{-1} \omega + {}^E \Omega^2 (M, \C [[\hbar ]])
$$
such that $d\theta = 0$
and let $\nabla_0$ be any ${\frak g}_0$-connection in $E$
There exists a
$\tilde{\frak g } $-valued $E$-form $A_{\theta}$ on $M$ such that 
$$
\nabla_{\theta} =\nabla_0 + A_{\theta},
$$
satisfying
$$
\nabla_{\theta} A_{\theta} +\frac{1}{2}[A_{\theta} ,A_{\theta} ] = \theta 
$$
and defining  a Fedosov connection on ${}^E${\bf W}. We will call $\theta
$ the {\bf curvature} of $\nabla_{\theta} $. The complexes
$$
({}^E \Omega (M, \vfA ) ,\nabla_{\theta} )
$$
and
$$
({}^E \Omega (M, {}^E Jets \otimes {}^E{\mathbb{W}}) ,\nabla_G
+\nabla_{\theta} ) 
$$
are acyclic in positive dimension,
$$
\App_E (M) \stackrel{\sim}{\rightarrow} Ker \left( \nabla_{\theta} |_{{}^E 
\Omega^0 (M, \vfA )}
\right) 
$$
is an $E$-deformation of the Poisson manifold $(M,\{ \ ,\ \}_E )$ with the
associated deformation of the jet bundle given by
$$
({}^E \underline{Jets} (M),*) \stackrel{\sim}{\rightarrow} Ker \left( (\nabla_G 
+\nabla_{\theta} )|_{{}^E
    \Omega^0 (M,  {}^E Jets \otimes {}^E{\mathbb{W}} )}
\right) . 
$$
\end{th}

\medskip

{\em Proof}.

\medskip 

1. {\bf Construction of Fedosov connection.}

\medskip

The construction of $\nabla_{\theta}$ is via recursion in the grading of
$\tilde{\frak g }$. Let  
$$
\nabla_{-1} = \nabla_0 + A_{-1} .
$$
Then 
$$
\nabla_{-1} ( A_{-1} ) \in \tilde{\frak g }_{-1} ,
$$
and hence 
$$
[A_{-1} , \nabla_{-1} ( A_{-1} )]=0.
$$
By the lemma \ref{acyclicity} above, there exists a $\tilde{\frak g
}_0$-valued one-form $A_0$ such that  
$$
\nabla_{-1} ( A_{-1} ) =[A_{-1} , A_0 ].
$$
Set
$$
\nabla_{00} = \nabla_{-1} + A_0 .
$$
We have
$$
[\nabla_{00} ,\nabla_{00} ]-\theta =0 \mbox{\ mod }\tilde{\frak g }_{\geq 0}
$$
Now suppose that we have constructed
$$
\nabla_{n} = \nabla_{-1} + A_{-1} + A_0 +\ldots +A_n ;\ A_p \in {}^E \Omega^1 (M, 
,
\tilde{\frak g }_p ) , 
$$
such that
$$
[\nabla_{n} ,\nabla_{n} ]-\theta =0 \mbox{\ mod }\tilde{\frak g }_{\geq n}.
$$
The ${}^E \Omega^3 (M, \tilde{\frak g }_{n-1})$-component of the identity
$$
[\nabla_{n},[\nabla_{n},\nabla_{n}] -\theta ]=0
$$
gives
$$
[A_{-1},([\nabla_{n},\nabla_{n}] -\theta)_{n} ]=0.
$$
Again by the lemma \ref{acyclicity} we can find an $A_{n+1}$ in ${}^E \Omega^1
(M , \tilde{\frak g }_{n+1})$ such that
$$
2[A_{-1},A_{n+1}] + ([\nabla_{n},\nabla_{n}] -\theta)_{n} = 0.
$$
But this means that, for $\nabla_{n+1}=\nabla_{n} + A_{n+1}$, 
$$
[\nabla_{n+1} ,\nabla_{n+1} ]-\theta =0 \mbox{\ mod }\tilde{\frak g }_{\geq
n+1}. 
$$
Since the filtration on $\tilde{\frak g }$ is complete, the above procedure gives $\nabla_{\theta}$ which obviously is a
Fedosov connection with curvature $\theta$.

\medskip

2. {\bf Acyclicity of $ ({}^E \Omega (M, {\mathbb A}^{\hbar }) ,\nabla_{\theta}
) $}.

\medskip

Let us write
$$
{}^E \Omega (M, {\hat{\mathbb A }}^{\hbar }) =\oplus_{k\geq 0} C^k
$$
where $C^k$ consists of $({\hat{\mathbb A }}^{\hbar })_k$-valued $E$-differential
forms on $M$. $\{ C^k \} $ is a complete filtration of our complex, and the
differential $\nabla_{\theta}$ on $C^k/C^{k-1}$ reduces to Ad($A_{-1}$). By
the lemma \ref{acyclicity} the corresponding spectral sequence degenerates and
hence the cohomology in positive dimensions is zero, while the kernel of
$\nabla_{\theta}$ is linearly isomorphic to $C^{\infty }(M)[[\hbar ]]$.

\medskip

3. {\bf Construction of the tensor $\varpi_E$}.

\medskip

Given a Fedosov connection $\nabla_{\theta }$ as constructed above, we get a
flat connection on the bundle of $\hat{\mathbb A}^{\hbar}$-valued $E$-jets:
\begin{equation}
\nabla : \Gamma ({}^E Jets \otimes {}^E{\mathbb{W}})\rightarrow {}^E \Omega^1
(M, {}^E Jets \otimes {}^E{\mathbb{W}})
\end{equation}
by setting
\begin{equation}
\nabla = \nabla_G \otimes 1 + 1 \otimes \nabla_{\theta }.
\end{equation}
Consider the embedding 
\begin{equation} \label{im}
({}^E \Omega ^*(M,  {}^E Jets), \nabla _G) \rightarrow ({}^E \Omega ^*(M,  {}^E Jets \otimes {}^E{\mathbb{W}}), \nabla _G \otimes 1 + 1 \otimes \nabla_{\theta})
\end{equation}
Note that this is a morphism of filtered complexes: the filtration on ${}^E Jets$ is by powers of the maximal ideal at any point, the filtrations on $\hat{\mathbb A}^{\hbar}$ and on $\wedge E^*$ are induced by their gradings, and the filtrations on the complexes in (\ref{im}) are tensor products of those filtrations.
Note that (\ref{im}) is a quasi-isomorphism because it induces a quasi-isomorphism of associated graded spaces of the above filtrations.

The fiberwise product gives us now an associative product $*$ on the space of
$E$-jets identified with the space of $\nabla$-flat sections:
$$
(l_1 *l_2) (D) = l_1 \otimes l_2 (\phi (D)),
$$
and hence by duality a map $\phi $ of left $C^{\infty } (M)$ modules (with $\C 
[[\hbar ]]$ as
the field of scalars)
$$
\phi : {}^E Op (M) \rightarrow  {}^E Op (M) \otimes_{C^{\infty } (M)}  {}^E Op
(M). 
$$
Since $\nabla$ commutes with $\nabla_G$, the Grothendieck connection acts on
the space of $\nabla$-flat sections as a derivation with respect to the
$*$-product, i.e. 
$$
\phi (XD )=(X\otimes 1 + 1 \otimes X ) \phi (D).
$$
But this means that the tensor
$$
\varpi_E = \phi (1) \in {}^E Op (M) \otimes_{C^{\infty } (M)}  {}^E Op (M)
$$
satisfies
$$
(l_1 *l_2) (D) = l_1 \otimes l_2 (\Delta_0 (D) \cdot \varpi_E ).
$$
It is now straightforward to see that the associativity of the $*$-product
implies that $\varpi_E $ satisfies the equation (\ref{assoc}). 

An immediate corollary of the proof above is the fact that any $E$-jet $l$ has
a unique continuation $\tilde{l}$ to a $\nabla$-flat section of the bundle of
$\W (E^*)$-valued jets. 

4. {\bf End of the proof}.

\medskip

To show that $\App_E (M)$ is an $E$ deformation of $\{ \ ,\ \}_E$, it is now
sufficient to prove that, for two $\C$-valued $E$-jets $l_1$ and $l_2$ 
$$
l_1 * l_2 = l_1 \cdot l_2 \mbox{mod}(\hbar ) \mbox{\ and }\\
(\frac{1}{\hbar }[l_1,l_2 ]) = \{l_1 ,l_2 \}_E  \mbox{mod}(\hbar ), 
$$
where ``$\cdot $'' is, as usual, the undeformed product and the commutator is
taken with respect to the deformed product. By the last sentence of the part
three of the proof, it is sufficient to show that, if $\tilde{l}_1$ and
$\tilde{l}_2$ are the extensions to flat $\nabla$ sections, then
$$
\tilde{l}_1 \tilde{l}_2 = (l_1 \cdot l_2 )^{~} \ \mbox{mod}(({\hat{\mathbb A
}}^{\hbar })_{\geq 1} )
$$
and
$$
(i\hbar )^{-1} [\tilde{l}_1 ,\tilde{l}_2] = \tilde{ \{ l_1 ,l_2 \}_E } \ 
\mbox{mod}(({\hat{\mathbb A }}^{\hbar })_{\geq 1} ).
$$
The first equation follows from the computation:
$$
\tilde{l}_1 \tilde{l}_2 =(l_1 +({\hat{\mathbb A }}^{\hbar })_{\geq 1})(l_2
+({\hat{\mathbb A }}^{\hbar })_{\geq 1}) =(l_1 \cdot l_2) +({\hat{\mathbb A 
}}^{\hbar
})_{\geq 1}. 
$$
To prove the second equation, we need a bit of notation. Let 
$$
e_1, \ldots e_n , f_1 , \ldots ,f_n 
$$
be a local symplectic basis of sections of $E$ and set, for any section $v$ of
$E$,
$$ 
v(m)^* =I_{{\omega}_m}(v(m)) \in {\hat{\mathbb A }}^{\hbar }(E_m).
$$
and 
$$
J: E_m \rightarrow E_m \ \mbox{with} \ J(e_i)=f_i,\ J(f_i) =-e_i
$$
We will write $\hat{v}$ for $v$ considered as a formal linear coordinate
function on $E_m^*$ and $\partial_{\hat{v}}$ as the fiberwise derivative on
the jet bundle. Modulo $({\hat{\mathbb A }}^{\hbar })_{\geq 1}$ we have
$$
\nabla_v = v + \partial_{J(v)\hat{} } + \partial_{v^*}
$$
and hence any $\nabla $-flat section $\tilde{l}$ satisfies
$$
\tilde{l}= l + \sum_{i=1, \ldots ,n} (\partial_{\hat{e}_i} l \cdot f_i^*
+\partial_{\hat{f}_i} l \cdot e_i^* ) \ \mbox{mod}(({\hat{\mathbb A }}^{\hbar 
})_{\geq
2}). 
$$ 
Since $[e_i^* ,f_j^* ] =\delta_{ij}$, we get, modulo $({\hat{\mathbb A }}^{\hbar
})_{\geq 1}$,
$$
\frac{1}{i\hbar }[\tilde{l}_1 ,\tilde{l}_2 ] \\
  = \sum_i (\partial_{\hat{e}_i} l_1 \ \partial_{\hat{f}_i} l_2
-\partial_{\hat{f}_i} l_1 \ \partial_{\hat{e}_i} l_2 ) \\
  = \{ l_1 ,l_2  \}_E
$$
This completes the proof of the theorem.

\bigskip

\section{Formal deformations associated to symplectic Lie algebroids, 
classification }
We will continue using notation from the previous section.

\subsection{Global structure of $E$-deformations}

Our next objective is to show that any $E$-deformation is in fact one of the
type constructed above and that the cohomology class of the curvature form of the associated Fedosov
connection is
a complete invariant of the deformation in our class. 

\bigskip

\begin{th}\label{global structure}
Let $\App (M)$ be an $E$-deformation for a symplectic Lie algebroid
$(E,[,],\rho ,\omega )$ on $M$, and $({}^E Jets  ,* )$ the associated deformed
bundle of algebras. There exists an isomorphism of bundles of algebras
$$
\phi : ({}^E Jets  ,* ) \rightarrow {}^E\mbox{\bf W}
$$
which maps the Grothendieck connection $\nabla_G$ to a Fedosov connection
$(\phi^{-1})^* (\nabla_G)$ on the Weyl bundle ${}^E \mbox{\bf W}$. The associated
deformation of the algebra of smooth functions 
$$
Ker\left( (\phi ^{-1})^* (\nabla_G)|_{\Gamma (M, {}^E\mbox{\bf W})} \right) 
$$
is isomorphic to $(\App (M),*)$.
\end{th}
{\em Proof}. 

We will begin by constructing the required isomorphism of bundles
locally. So let $U$ be an open subset of $M$ on which $E$ admits a symplectic
basis 
$$
(e_1 ,f_1 ,\ldots ,e_n ,f_n ) .
$$
In this basis $E|_U$ becomes identified with 
$$
U\times \R^{2n}
$$
and we denote by $\R^{2n}$ the linear subspace of sections of $E$ of the form
$$
m\mapsto (m, v), \ v \ \mbox{fixed in $\R^{2n}$} .
$$
For a section $v$ of $E$ over $M$ we will denote by $l_v$ the
$E$-jet given by  
$$
\begin{array}{l}
l_v (1)=0 ,\ l_v (w)=\omega (v,w),\\
l_v (Dw) = 0 \ \mbox{whenever $D\in {}^E Op^+$ }.
\end{array}
$$
Clearly the set 
$$
\{ l_v \mid v\in \R^{2n} \ \}
$$
generates the algebra of $E$-jets modulo $\hbar$.

Note that
$$
[l_v ,l_w] (1)=i\hbar \omega (v,w) -i \hbar^2 B_2(v,w) +\ldots
$$
where $B_2$ is a skew-symmetric bilinear complex-valued form on
$\R^{2n}$. Since $\omega $ is non-degenerate, there exists a linear
transformation $A$ of $\R^{2n}\otimes_{\R} \C$ such that
$$
-B_2 (v,w) = \omega (Av ,w) +\omega (v,Aw ).
$$
Denoting by $\tilde{v}$ the section
$$
\tilde{v} =v +\hbar A(v) 
$$
we get
$$
[l_{\tilde{v}} ,l_{\tilde{w}}] (1)=i\hbar \omega (v,w) -i \hbar^3  B_3(v,w)
+\ldots 
$$
An obvious induction gives now an invertible map
\begin{equation}
\begin{array}{ccc}
\R^{2n}\otimes_{\R} \C [[\hbar ]] & \rightarrow  & \R^{2n}\otimes_{\R} \C [[\hbar 
]] \\
v  \otimes 1     &  \mapsto       &\tilde{v} =v\otimes 1+O(\hbar )
\end{array}
\end{equation}
such that
\begin{equation}
[l_{\tilde{v}} ,l_{\tilde{w}}] (1)=\hbar \omega (v,w).
\end{equation}
Thus we get isomorphisms 
$$
\phi_i \;: \; {}^E Jets | U_i \stackrel{\sim}{\rightarrow} {}^E {\mathbb{W}} | U_i
$$
for a good cover. The transition isomorphisms
$$
g_{ij} = \phi_i^{-1}  \phi_j \in \check{C}^1 (M,Aut^1 ({}^E \mbox{\bf W}))
$$
take values in the pronilpotent group $G^{\geq 1}$ of automorphisms of ${\mathbb{W}}$ which preserve the filtration and are equal to the identity on the associated graded space. Therefore there is a global isomorphism of filtered algebras 
$$
{}^E Jets  \stackrel{\sim}{\rightarrow} {}^E {\mathbb{W}}
$$
The image of $\nabla _G$ under this isomorphism is a Fedosov connection.

{\bf The characteristic class }

Consider the Lie algebra central
extension: 
\begin{equation}
0\rightarrow (i\hbar )^{-1} \R +\C [[\hbar ]] \rightarrow \tilde{\frak g} 
\rightarrow {\frak g }
\rightarrow 0. 
\end{equation}
For a Fedosov connection $\nabla$, let $\tilde{\nabla}$ be any lifting of $\nabla$ to a $\tilde{\frak{g}}$-valued connection. Then
$$ \theta = \frac{1}{2} [\tilde{\nabla},\tilde{\nabla}]$$
 is an element of ${}^E \Omega ^2 (M, (i\hbar )^{-1} \R +\C [[\hbar ]])$ such that $d \theta = 0$.
\begin{definition}
\begin{em}
The characteristic class of the deformation $\App (M) $ is the cohomology class
$$
\theta \in{}^E H^2(M,(i\hbar )^{-1} \R +\C [[\hbar ]] ). 
$$
\end{em}
\end{definition}


\subsection{Classification of $E$-deformations}

\begin{th}
Let $(E,[\ ,\ ],\rho ,\omega ) $ be a symplectic Lie algebroid on $M$.

The characteristic class of an $E$-deformation is well defined.
Two $E$-deformations $\App_1 (M)$ and $\App_2(M)$ are
isomorphic if and only if their characteristic classes $\theta_1$ and
$\theta_2$ are equal.  Thus, the affine space 
$$
\frac{1}{i\hbar}\omega + {}^E H^2 (M,\C [[\hbar ]]) 
$$
completely classifies $E$-deformations of $M$ up to isomorphism.
\end{th}
{\em Proof}. 

By the theorem \ref{global structure} we can assume that both
deformations are given by Fedosov construction with connections $\nabla_1$ and
$\nabla_2$ on the Weyl bundle ${}^E \mbox{\bf W}$. Let  $\tilde{\nabla_1}$ and
$\tilde{\nabla_2}$ be their liftings. Note that the characteristic classes are
given by the curvatures 
$$
\theta_i =\frac{1}{2}[\tilde{ \nabla}_i ,\tilde{\nabla}_i ], \ i=1,2,
$$
{\bf 1}. We assume that these characteristic classes are cohomologous. 

 Let
$$
\theta_1 -\theta_2 ={}^E d\alpha , \alpha \in C^{\infty } (M)[[\hbar ]] ).
$$
But then, replacing $\nabla_2$ by $\nabla_2 +\alpha $, we get two connections
with the same curvature and unchanged deformations. So we can also assume that
$\theta_1 =\theta_2$ at the level of forms. We will construct an element of $Aut^1 ({}^E \mbox{\bf
W})$ which conjugates the two connections.

So, let
$$
\tilde{\nabla_1} = A_{-1} + \nabla_0 + A_1 + \ldots
$$
where $A_i$ are $\tilde{\frak{g}}_i$-valued one-forms on $M$ and $\nabla_0$ is induced by a unitary connection in $E$ (note that $\tilde{\frak{g}}_0 = {\frak{g}}_0 \oplus {\mathbb{C}}$ canonically). Let 
$$
\tilde{\nabla_1} - \tilde{\nabla_2} = R_0 +R_1 + \ldots ,\ R_i \in  {}^{E} \Omega^1  (M,
\tilde{\frak g}_i). 
$$
The equality of the curvatures of the two connections gives
$$
[A_{-1}, \nabla_0 +R_0 ]=0.
$$
Since the ad$A_{-1}$-complex is contractible by Lemma \ref{acyclicity}, we can
find $\delta_1$ in ${}^E \Omega^0 (M,\tilde{\frak g}_1 )$ such that 
$$
R_0 = [\delta_1 ,A_{-1} ].
$$
Replacing $\nabla_2 $ by Ad(exp$\delta_1$)($\nabla_2$), we get
$$
\nabla_1 - \nabla_2 =0 \ \mbox{mod} \tilde{\frak g}_{\leq 0}.
$$
Continuing in this vein, the induction on the grading of $\tilde{\frak g}$
gives us a sequence 
$$
\delta_i \in {}^E \Omega^0 (M,\tilde{\frak g}_i ) , \ i=1,2,\ldots ,
$$
such that
$$
\ldots
\mbox{Ad(exp$\delta_3$ )}\mbox{Ad(exp$\delta_2$)}\mbox{Ad(exp$\delta_1$)}
(\nabla_2 ) = \nabla_1 .
$$
Hausdorff-Campbell formula implies now that there exists an element of $Aut^1
({}^E \mbox{\bf W})$ conjugating the two connections, and hence the two
$E$-deformations are isomorphic.

{\bf  2}. Suppose now that the two deformations are isomorphic. 

This means that there exists an isomorphism of the deformed jet bundles: 
$$
({}^E Jets ,*_1 )\stackrel{\Phi }{\rightarrow} ({}^E Jets ,*_2 ).
$$
such that
$$
\Phi =id + O(\hbar ).
$$
But this implies that the curvature forms of the two corresponding
connections are cohomologous and we can apply the result above. 

This finishes the proof of the theorem. 

\medskip

The following is an immediate corollary of the previous two sections.

\medskip

{\bf Structure of derivations}

\medskip

\begin{th}\label{structure of derivations}
Let $\App (M)$ be an $E$-deformation of $M$ given by a Fedosov connection
$\nabla$ on the Weyl bundle ${}^E${\bf W}. Any $E$-derivation of $\App (M)$
extends to a derivation of the Weyl bundle which maps $\nabla$-flat sections
to $\nabla$-flat sections. In particular, the space ${}^E Der (\App (M) )$ of
E-derivations of $\App (M)$ is in bijective correspondence with
$$
\{ \ l\in \Gamma (M,{}^E \mbox{\bf W})\ | \ \nabla (l) \ \mbox{is center
valued}\ \} . 
$$
There exists a Lie algebra extension
$$
0\rightarrow \mbox{Ad}(\App (M) )\rightarrow {}^E Der (\App (M) )\rightarrow
{}^E H^1 (M,\C [[\hbar ]] )\rightarrow 0.
$$
\end{th}
{\em Proof}. All the statements above follow immediately from the fact that
any $E$-derivation of the $E$-deformation $\App (M)$ extends by definition to a
derivation of the associated deformed $E$-jet bundle. 

\medskip


\subsection{Gelfand-Fuks construction}\label{Gelfand-Fuks construction}


Suppose that $(E,[\ ,\ ],\rho ,\omega $ is a symplectic Lie algebroid and that
$(\App (M),*)$ is an $E$-deformation of $M$. Let  ${}^E${\bf W} be a Weyl
bundle and $\nabla$ a Fedosov connection associated to this deformation.
 Choose any local trivialization of the bundle $E$ on any open subset $U$ of $M$.  Let the Fedosov connection be of the form ${}^Ed + A_U$ in this trivialization. The flatness of $\nabla$ translates into
\begin{equation}\label{flat}
{}^{{E}}  dA_{U} +\frac{1}{2} [A_{U}, A_{U} ]=0
\end{equation}
and in particular implies that
\begin{equation}
({}^E \Omega (M,{\mathbb L }), {\nabla} )
\end{equation}
is a complex. In the future we will use $\nabla_{\mathbb L} $ to denote
$\nabla$ acting on this complex.

Recall that for any Lie algebra ${\frak{g}}$, a Lie subalgebra ${\frak {h}}$, and any ${\frak{g}}$ module $L$ the complex $(C^*({\frak g},{\frak h};{ L })$ of relative Lie algebra cochains is defined.

\medskip

\begin{definition} \label{gf}
\begin{em}
Let ${\mathbb L}$ be a continuous $({\frak g },\mbox{U(n)} )$-module. The
{\bf Gelfand-Fuchs map} is the map of complexes: 
$$
gf:(C^*({\frak g},{\frak u}(n);{\mathbb L }), \partial_{Lie} ) \rightarrow ({}^E
\Omega (M,{\mathbb L }) ,\nabla_{\mathbb L} ) 
$$
which is defined as follows. Let $l$ be a $k$-cochain of the relative Lie algebra
complex, and let $e_1,\ldots ,e_k $ be sections of $E$. We set
$$
gf(l)(e_1,\ldots ,e_k)=l(A_U({e}_1 ),\ldots ,A_U({e}_k )).
$$
Note that since $l$ is a relative cochain, the result is independent of the
choice of the trivialization and that the equation (\ref{flat}) implies that 
$$
gf \circ \partial_{Lie}  =\nabla_{\mathbb L} \circ gf .
$$
\end{em}
\end{definition}



\subsection{Example: symplectic manifolds}

Let $(M,\omega)$ be a symplectic manifold, and $\{ \ ,\ \}$ the associated
Poisson bracket on $M$. For simplicity we will assume through the rest of this
section that $M$ is connected. As the symplectic Lie algebroid we will take
the sheaf of all vector fields on $M$. The results of the previous sections
can be  formulated as follows (cf. \cite{LW}, \cite{19 Fedosov's book}, \cite{De 2}, \cite{NTfamilies}).
   \begin{th}
The set of isomorphism classes of formal deformations of $C^{\infty}(M)$ with
$$
f*g =fg + O(\hbar ), \ f,g\in C^{\infty}(M)
$$
and
$$
[f,g]=i\hbar \{ f,g \} +O(\hbar^2 ), \ f,g\in C^{\infty}(M)
$$
is in bijective correspondence with the elements $\theta$ of the 
space 
$$
(i\hbar )^{-1} \omega + H^2 (M,\C [[\hbar ]]).
$$
Every such deformed algebra is isomorphic to 
$$
Ker \nabla_{\theta} |_{\mbox{\bf W}}
$$
for a Fedosov connection on the Weyl bundle {\bf W}.
   \end{th}
\medskip

The structural results from the previous section give us the following
corollaries.
\begin{corollary}
Let $\App (M)$ be a formal deformation of a symplectic manifold $(M,\omega
)$. There exists a unique up to a scalar multiple $\C [\hbar^{-1} ,\hbar
]]$-valued trace on $\App_c (M)$. Up to normalization factor this
trace has the form
$$
Tr(f) =\frac{1}{(\frac{1}{2} dimM)! (i\hbar)^{\frac{1}{2} dimM}}\int_M f
\omega^{\frac{1}{2} dimM} +O(\hbar^{-\frac{1}{2} dimM+1})
$$
\end{corollary}
{\em Proof.} Recall that a deformation of a symplectic manifold is locally
unique (the characteristic class lies in $H^2 (U)$ which vanishes for each open
contractible subset $U$). Let $\{ U_i \}_{i\in I}$ be a locally finite
covering of $M$ by open contractible subsets and $\{ \rho_i \}$ an associated
partition of unity. By the theorem \ref{local trace}, we get the family of
traces $Tr_i$ on $\App_c (U_i)$. Set  
$$
Tr (f)=\sum_i Tr_i (\rho_i f).
$$
This is a well-defined trace (cf. \cite{19 Fedosov's book}, \cite{NT1} for the proof).
\medskip

For completeness let us record the following result which is a stronger version of 
[25]: here by isomorphism we mean an isomorphism of unital algebras.

\begin{th}
Given a formal deformation $\App (M)$ of a compact symplectic manifold
$(M,\omega )$, there exists within the isomorphism class of $\App (M)$ a
$*$-product on $C^{\infty }(M)[[\hbar ]]$ such that the above trace has the
form 
$$
Tr(f) =const \cdot \frac{1}{(\frac{1}{2} dimM)! (i\hbar)^{\frac{1}{2}
    dimM}}\int_M f \omega^{\frac{1}{2} dimM} .
$$
In the terminology of \cite{CFS} this is a closed deformation.
\end{th}
{\em Proof.} The canonical trace constructed above has the form
$$
f\rightarrow c(\hbar )\int_M^{} T(f)\omega^n ,
$$
where $n$ is half the dimension of $M$ and $T$ is a linear transformation acting
on the space of smooth functions and of the form
$$
T(f)=fT(1), \ T(1)=1 + O(\hbar )\in C^{\infty }(M)[[\hbar ]].
$$
What we need to find is a linear transformation $S$ of smooth functions on $M$
such that 
\begin{itemize}
\item $f\rightarrow \int_M^{}  (T+\hbar S)(f)\omega^n $ is a trace with
  respect to the original $*$-product
\item $\hbar S(1)=1 -T(1) $
\end{itemize}
Once this is done, the new product will be given by
$$
f*_{new} g =(T+\hbar S)((T+\hbar S)^{-1}(f)*(T+\hbar S)^{-1}(g)).
$$ 
By the uniqueness of the trace, the first condition above is equivalent to  
$$
\int_M (T+\hbar S)(f)\omega^n =k \int_M (T)(f)\omega^n, \ k=1+O(\hbar )
$$
i.e.
$$
(T^* + \hbar S^*)\omega^n  = k T^* \omega^n 
$$
where the adjoint ${}^*$ is taken with respect to the duality between
smooth functions and 2n-forms given by integration:
$$
C^{\infty} (M)\times \Omega^{2n}(M) \ni (f,\nu  ) \mapsto \int_M f\nu .
$$
In particular, we need to solve the equations
$$
S^* \omega^n = \frac{k-1}{\hbar} T^* \omega^n =\frac{k-1}{\hbar} T(1),
$$
subject to the condition
$$
S(1)=\frac{1-T(1) }{\hbar } .
$$
We will ask for $S$ of the form 
$$
S=\frac{1-T(1) }{\hbar } + X, \ X\in \Gamma (M,TM).
$$
Since
$$
X^* \omega^n =L_X \omega^n ,
$$
the above equation for $S^*$ reduces to
$$
L_X \omega^n = \frac{kT(1)-1}{\hbar } \omega^n .
$$
Since the map
$$
Vect(M)\ni X \rightarrow i_X \omega^n \in \Omega^{2n-1} (M)
$$
is an isomorphism, a 2n-form on $M$ is exact precisely when it is of the form
$L_X$ for some vector field $X$ on $M$. In particular our equation has a
solution for the constant $k$ given by 
$$
k^{-1} [\omega^n ] =[T(1)\omega^n ]\in H^{2n}(M,\C [[\hbar ]]),
$$
which is in $1+\hbar \C [[\hbar ]]$ since $T(1)$ is a formal power series
with leading coefficient 1.

The $*_{new}$ associated to this solution satisfies the claim of the theorem.

\section{Complex symplectic manifolds}  \label{complex symplectic manifolds}

Let $M$ be a complex manifold. We will denote by ${\cal O}_M$
the structure sheaf of holomorphic functions on $M$ and by $\cal O_{\infty }$
the sheaf of smooth functions on $M$.
\begin{definition}
A {\em deformation quantization} of a manifold $M$ is a formal one parameter
deformation of the structure sheaf ${\cal O}_M$, i.e. a sheaf of algebras
${\mathbb A}^\hbar_M$ flat over $\mathbb C [[\hbar]]$ together with an isomorphism 
of sheaves of
algebras $\psi: {\mathbb A}^\hbar_M\otimes_{\mathbb C [[\hbar]]}{\mathbb C }
\rightarrow {\cal O}_M$.

The formula
\[
\lbrace f, g\rbrace = \frac{1}{\hbar}[\tilde f,\tilde g ]
+ \hbar\cdot {\mathbb A}^\hbar_M\ ,
\]
where $f$ and $g$ are two local sections of ${\cal O}_M$ and $\tilde f$, $\tilde 
g$
are their respective lifts to ${\mathbb A}^\hbar_M$, defines a Poisson structure on $M$
called the Poisson structure associated to the deformation quantization
${\mathbb A}^\hbar_M$.

The deformation quantization ${\mathbb A}^\hbar_M$ is called {\em symplectic} if the 
associated Poisson structure is nondegenerate. In this case $M$ is symplectic, i.e., has a holomorphic symplectic form.
In what follows we will only consider symplectic deformation quantizations,
so assume that ${\mathbb A}^\hbar_M$ is symplectic, and $\omega $ denotes the
associated symplectic structure on $M$. 
\end{definition}

Let us note first that, given a deformation ${\mathbb A}^\hbar_M$ as above, $\psi$
induces locally an isomorphism of sheaves of $\C$-vector spaces:
$$
\tilde{\psi }: \prod (\hbar^{n+1} {\mathbb A}^\hbar_U )
/\hbar^{n} {\mathbb A}^\hbar_U \simeq {\cal O}_U [[ \hbar ]].
$$
Since locally there is no cohomology, this implies that there exist local
isomorphisms
$$
\Phi_U :{\mathbb A}^\hbar_U \rightarrow  ({\cal O}_U [[ \hbar ]] ,*_{U} )
$$
of $\k$-algebras with transition isomorphisms $G_{UV} = \Phi_U \Phi_V^{-1}$ of
the form
\begin{equation}
G_{UV} = id +\hbar D_1^{UV} +\hbar^2 D_2^{UV} +\ldots .  
\end{equation}

In the rest of this section we will work under following continuity condition.
\begin{assumption}
Both the local products $*_U$ and the linear transformations $D_i^{UV}$ above
are given by holomorphic (bi-)differential operators. 
\end{assumption}

\subsection{Weyl bundle, Fedosov connection, etc.}

Let us start with a few constructions associated with complex symplectic
manifolds which are analogous to the smooth case.

Let $\fA$ denote the Weyl algebra
over $\C [[\hbar ]]$ of the standard symplectic structure $(\C^{2n},
\omega_{st})$, {\bf and set, just for this section, }
\begin{equation}
  \label{complex3}
  \tilde{\frak g}  =\{ \hbar^{-1}f | f\in \fA\}
\end{equation}
with the Lie bracket given by the commutator in $\fA$. We will denote by $\frak
g$ the quotient Lie algebra $Der \fA$. Note that ${\frak {sp}} (2n, \C )$ is a
subalgebra of $\tilde{\frak g}$ and that its adjoint action integrates to the
action of G=Sp(2n, $\C$) on $\tilde{\frak g}$. 

Let $\cal P$ denote the principal Sp(2n,$\C$)-bundle of symplectic frames in
the holomorphic tangent bundle $T=T^{1,0} (M)$, with the complex structure
induced in the obvious way from the complex structure on the complex Lie group
G. We define $\vfA$-valued differential forms by
\begin{equation}
  \label{complex4}
  \Omega^{p,q} (M, \vfA )=(\Omega^{p,q}({\cal P}) \otimes \vfA )^{basic}.
\end{equation}
Note that $\overline{\partial}$ extends automatically to give an analogue of
the Dolbeault complex
$$
(\Omega^{p,*}(M,\vfA ), \dbar )
$$
for all p.

We denote by {\bf W} the (holomorphic) bundle
${\cal P} \times_{Sp(2n,\C)} \vfA$ and by ${\cal G}$ the gauge group of fiberwise inner
automorphisms of {\bf W}.  

Let $\nabla_0$ be any Sp(2n,$\C$)-connection of type (1,0) in {\bf W}; locally 
$$
\nabla_0 = \partial +\mbox{ad} \alpha \ :\Omega^{p,q} \rightarrow \Omega^{p+1,q}
$$
with $\alpha $ a ${\frak  sp} (2n, \C)$-valued form of type (1,0)
on $M$.

Let $A_{-1}$ denote the canonical ${\bf W}_{-1}$-valued holomorphic one-form on $M$: 
$$A_{-1} : T \stackrel{\omega}{\rightarrow} T^* \stackrel{\sim}{\rightarrow} {\bf W}_{-1}$$

\begin{definition}
\begin{em}
A Fedosov connection is a connection on {\bf W } of the form
$$
\nabla = \dbar + \nabla_0 + ad A + ad B,
$$
where $\nabla_0$ is a connection on $TM$ preserving both the complex and
symplectic structure, hence in a local symplectic frame of the form
\begin{equation}
  \label{conn}
\partial + \dbar +\alpha_U  
\end{equation}
 with $\alpha_U \in \Omega^{1,0} (U, {\frak sp}\mbox{(n)})$,
$$
A = A_{-1} + A_0 + \ldots ,\ A_i \in \Omega^{1,0}(M,\g_i )
$$
and
$$
B =B_1 + B_2 +\ldots ,\ B_i \in \Omega^{0,1}(M,\g_i ),
$$
which satisfies the flatness condition
$$
\nabla^2 = 0
$$
\end{em}
\end{definition}

Note that the curvature $\nabla^2$ of a Fedosov connection splits into
the components 

  \label{complex9}
  \begin{enumerate}
  \item $\nabla_0 A +\frac{1}{2}[A,A] =\theta \in \Omega^{2,0}(M,\k ); $
  \item $\dbar A + \nabla_0 B +[A,B]  =\phi \in \Omega^{1,1}(M; \k); $
  \item $\dbar B + \frac{1}{2}[B,B] =\tau \in \Omega^{0,2} (M, \k ),$
   \end{enumerate}
and it satisfies the Bianchi identity
 $$
 d (\theta + \phi + \tau) = [\nabla ,\nabla^2 ]=0 .
 $$

\begin{th}[Classification of Fedosov connections]
  \label{ccc5}
Let $(M,\omega )$ be a complex symplectic manifold such that the inclusion of
sheaves $\C_M \rightarrow {\cal O}_M$ induces a surjection
$$
H^1 (M,\C ) \rightarrow H^1 (M, {\cal O}).
$$
Let $\nabla$ and $\nabla^{'}$ be two $\tilde{\frak g}$-valued Fedosov
connections on the associated Weyl bundle {\bf  W}. Then $\nabla$ and
$\nabla^{'}$ have the same curvature class in $H^2 (M,\C )$ if and only if
there exists a $\k$-valued one-form $\alpha$ such that $\nabla$ and
$\nabla^{'} +\alpha$ are conjugate by a gauge transformation by an element of
$ \Gamma (M, \exp {\frak g}_{\geq 1})$. 
\end{th}

{\it Proof}. Suppose first that $\nabla^2 = (\nabla^{'})^2$. We will work by 
induction on $n$, where 
\begin{equation}
  \label{ccc1}
  \begin{array}{rl}

  \nabla_n = A_{-1}+\nabla_0 + A_1 & + \ldots +A_n + \\
                      \dbar +B_1 & + \ldots +B_{n} + B_{n+1}
\end{array}
\end{equation}
So suppose that $\nabla_n =\nabla^{'}_n +\sum_{i\leq n}\alpha_i$, where 
$\alpha_n$ is the component of $\alpha$ in 
$\Omega^{10}(M, \fA_n ) \oplus \Omega^{01}(M, \fA_{ n+1} )$. 
The fact that the curvature forms coincide implies the identities 
\begin{itemize}
\item $[A_{-1} , A_{n+1}] = [A_{-1} , A_{n+1}^{'}]  + \partial \alpha^{10}_n$,
\item $[A_{-1} , B_{n+2}]-\dbar ( A_{n+1}) = 
 [A_{-1} , B_{n+2}^{'}] -\dbar (  A_{n+1}^{'})  + \partial \alpha_{n}^{01}$, 
\item $\dbar (B_{n+2}- B_{n+2}^{'})=0$.
\end{itemize}

Suppose first that n is even. Since ad$A_{-1}$ is acyclic, we can find an
$x_{n+2} \in \Gamma (M,\fA_{n+2})$ such that $ A_{n+1} - A_{n+1}^{'} = [A_{-1}
, x_{n+2}]$. But then the gauge transformation Ad$e^{x_{n+2}}$ allows us to set $
A_{n+1} = A_{n+1}^{'}$. Now the second equation implies that $ B_{n+2}$
and $B_{n+2}^{'}$ differ by a scalar-valued form $\beta$ which satisfies
$\dbar \beta =0$. By our assumption, there exists a $\hbar^\frac{n+2}{2}
\C$-valued section $y_{n+2}$ and a $d$-closed scalar-valued one-form
$\alpha_{n+1}$ such that $\beta_{n+2} = \dbar y_{n+2}
+\alpha_{n+1}^{01}$. In particular, using gauge transformation  Ad$
(exp(-y_{n+2}))$ we can assure that 
$$
 (\nabla^{'})_{n+1}= \nabla_{n+1} + \alpha_{n+1}^{01} .
$$ 

Suppose now that n is odd.

Since ad$A_{-1}$ is acyclic on $\Omega^{10}(M,\fA
)$, the first identity implies that there exists a $\fA_{n+2}$-valued section
$x_{n+1}  $ such that
$$
A_{n+1}+\alpha_{n+1}^{10} - A_{n+1}^{'} = [A_{-1} , x_{n+2}]
$$
Using gauge transformation $ Ad (exp(x_{n+2}))$ we can replace $\nabla^{'}$ by
a connection for which the $A$-components agree up to the order $n+1$ up to the
appropriate component of $\alpha$. But then the second identity implies that
that $ B_{n+2}$ and $B_{n+2}^{'}$ differ by a scalar-valued one-form of odd
degree and hence coincide.

Since in the curvature of a Fedosov connection is gauge invariant, we proved
that when the curvature forms coincide, the two connections are gauge
equivalent up to a closed scalar-valued one-form. Since changing the
connection by a scalar one-form changes the total curvature by the differential of this form, the result follows. 

\begin{remark}  \label{gauge partial}
\begin{em} $\mbox{}$
A more precise statement given by the above proof is the
following. Let 
$$ \theta_{(n)} = \sum _{i\leq n-1} \theta_i + \sum _{i\leq n}\phi_i + \sum_{i\leq n+1} \tau_i$$
Suppose that we are given two Fedosov connections $\nabla$ and
    $\nabla^{'}$. If their curvatures $\theta$, $\theta '$ satisfy $\theta_{(n)} = \theta_{(n)}'$ then there exists a Fedosov connection $\nabla ''$ which is gauge equivalent to  $\nabla '$  such that $\nabla '_ {(n)} = \nabla ''_{(n)}$.
\end{em}
\end{remark}

\subsection{The structure of formal deformations}

\medskip

\begin{th}  \label{structure of formal deformations}
Given a symplectic deformation of $M$ there exists a Fedosov
connection 
  $\nabla$ on the Weyl bundle {\bf W} such that the sheaf ${\mathbb A}^\hbar_M$
is isomorphic to the sheaf
$$
U \rightarrow  Ker\nabla|_{\Omega^0 (U,\fA)}
$$ 
\end{th}
{\it Proof.} 
As a direct consequence of our assumption, the local product $*_U$ extends to
give a graded algebra structure to the Dolbault complex
$$
{\cal O}_U [[ \hbar ]] \hookrightarrow (\Omega^{0,*} (U)[[\hbar ]],\dbar , *_U
),
$$
where functions of $\overline{z}$ only (i.e. antiholomorphic) and the
differentials $d\overline{z}$ are central and the $\dbar$ operator
acts as an odd derivation with square zero. Moreover both product and $\dbar$
commute with the action of the transition functions $G_{UV}$ and so these local
complexes glue together to give a resolution of the sheaf of algebras $\App_M$
of the form 
\begin{equation}
 \App_M \hookrightarrow ({\cal F}^* ,\dbar ).
\end{equation}

The sheaves ${\cal F}^i$ are locally isomorphic and hence isomorphic to
the fine sheaf $\Omega^{0,i}(M,\k )$. Using this isomorphism we get on
$\Omega^{0,*}(M,\k )$ the following structures.
\begin{enumerate}
\item A structure of a graded algebra with an associative product $*$ given by
bidifferential operators constructed out of the vector fields of type (1,0)
(i.e. a $T^{10}$-deformation).
\item An odd derivation (with respect to the deformed product $*$)
  $\dbar^{\hbar}$ satisfying 
$$
\dbar^{\hbar} =\dbar + 0( \hbar ) \ \mbox{ and }\ (\dbar^{\hbar})^2 =0.
$$
\end{enumerate}

Once we have the ${\cal T}^{1,0}$-deformation of the sheaf
$(\Omega^{0,*} (M)[[\hbar ]],\dbar^{\hbar})$, the construction of the
jet-bundle, the associated Grothendieck connection $\nabla_G$ and the
identification of the jet bundle with the Weyl bundle $\Omega^{**}(M,{\bf W})$
associated to $\omega$ goes through word for word as in the smooth case. Since
$\dbar^{\hbar}$ acts as a derivation, it lifts to a derivation $\nabla^{01}$
of the jet bundle which commutes with the Grothendieck connection and has
square zero. In particular, the image of $\nabla_G + \nabla^{01}$ under the isomorphism of the jet bundle with the Weyl
bundle defines a Fedosov connection on $\Omega^{**}(M,{\bf W})$ with required
properties. 

\medskip

For the ease of the reader we will sketch a more explicit construction of the
Fedosov connection below.

(i) Local jet bundles.

\noindent For a local coordinate neighbourhood $(U,z_1 ,\ldots , z_{2n})$ we identify
$$
Jets_U = U \times \k [[\fz_1 ,\ldots  ,\fz_{2n} ]].
$$
We will denote (holomorphic) sections of the jet bundle by functions $f(z,\fz
)$, i.e. a formal power series in the (commuting) formal variables
$\fz_i$ with coefficients holomorphic functions in $z_i$. For any holomorphic function $F$ we put
$$
F(z+\fz )=\sum \frac{\fz^{\alpha} }{\alpha !} \partial_z ^{\alpha } F (z).
$$
For any holomorphic differential operator $D = \sum P_{\alpha}(z) \partial ^{\alpha}_z$, put
$$
\widehat{D} = \sum P_{\alpha}(z + \fz) \partial ^{\alpha}_{\fz}
$$
The 
section $\fz_i$ corresponds to the functional on holomorphic differential
operators given by
$$
Op_{hol} (U) \ni D \mapsto \widehat{D}\fz_i | _{\fz=0} \in {\cal O}_U .
$$
The $*_U$ product is according to our assumption given by an expression of the
form 
$$
f*_U g = \sum F_{\alpha ,\beta } (z,\hbar ) \partial_z^{\alpha} f
\partial_z^{\beta} g 
$$
and we set, for the jets $\phi ,\psi $ on U,
$$
\phi *_U \psi = \sum F_{\alpha ,\beta } (z+\fz ,\hbar )
\partial_{\fz}^{\alpha} \phi \partial_{\fz}^{\beta} \psi .    
$$
 
The Grothendieck connection is in our local coordinates given by the
expression
$$
\nabla_G f = \sum_i dz_i (\partial_{z_i} - \partial_{\fz_i} )f .
$$

(ii) Global jet bundle

\noindent Given a transformation $G: {\cal O}_U [[\hbar ]] \rightarrow {\cal O}_U
[[\hbar ]]$ of the form
$$
(Gf)(z) = \sum T_{\alpha }(z, \hbar ) \partial_z^{\alpha} f(z),
$$
we define its jet by
$$
jet (G) = \sum T_{\alpha }(z+\fz , \hbar ) \partial_{\fz}^{\alpha}
$$
Let now, for a point $m$ in $U$,
$$
\phi_U : \k [[\fz ]] \rightarrow Jets_m
$$
be the isomorphism provided by the coordinates $\{ z_i^U \}$ in $U$. Then we
glue the local jet bundles using the transition functions
$$
\phi_U^{-1} jet (G_{UV}) \phi_V : Jets_U |_{U\cap V} \rightarrow  Jets_V
|_{U\cap V}. 
$$
It is immediate from the construction that we get a bundle
$Jets (\A )$ of algebras on $M$ (i.e. the local products $*_U$ define a
global product), carrying a flat connection $\nabla_G$ of the form $\partial +
$derivation. Since both the transition functions and the product are given
locally by holomorphic differential operators, we get the associated complex
of sheaves of algebras.
$$
( \Omega^{*,*} (M ,Jets \A ) ,\nabla_G +\dbar ).
$$

(iii) The Fedosov connection

\noindent Locally the definition of Weyl bundle gives isomorphisms
$$
Jets_U \stackrel{\sigma}{\rightarrow}  {\bf W }_U
$$
with $\sigma \in \Gamma_{hol} (U, \mbox{exp} {\frak g}_{\geq 1} )$. Using
completeness of {\bf W} in the filtration of $ \tilde{\frak g}$, these local
isomorphisms give rise to an isomorphism of the associated smooth vector
bundles
$$
\Omega^{*,*} (\cdot ,Jets(\A ) )\rightarrow \Omega^{*,*} (\cdot , \vfA ) .
$$
Under this isomorphism $\nabla_G + \dbar$ gives a Fedosov connection
$\nabla_F$ such that 
$$
\A_M \simeq \mbox{Ker}\nabla_F |_{\Omega^{0,0}} .
$$

\bigskip
As a corollary we get the following result.
\begin{th} 
Let $(M,\omega )$ be a complex symplectic manifold for which the map $H^1(M,\C
) \rightarrow H^1 (M; {\cal O}_M)$ is surjective. Two formal
deformations of $(M,\omega )$ with the same cohomology class of the curvature
of the associated Fedosov connection are isomorphic.   
\end{th}
{\it Proof}. Since by above a formal deformation of $(M,\omega )$ is
automatically of the form
$$
U \rightarrow  Ker\nabla|_{\Omega^0 (U,\fA)}
$$ 
for some Fedosov connection, the result follows from theorem \ref{ccc5}

\begin{corollary}[local structure of deformations]\label{complex18}
Any formal deformation of a complex symplectic structure is locally isomorphic
to the sheaf of holomorphic functions on an open subset of $\C^{2n}$ endowed
with Weyl product.
\end{corollary}
{\it Proof}. Since any formal deformation comes from a Fedosov connection and
is uniquely determined by its curvature class, it is locally isomorphic to the
deformation of $\C^{2n}$ with its standard symplectic structure and with
respect to any Fedosov connection. Let complex coordinates in
$\C^{2n}$, $(z_1 , \eta_1 , \ldots , z_n ,\eta_n )$ be such that 
$$
\omega = \sum dz_i \wedge d \eta_i 
$$
The expression
$$
\nabla = d - (i\hbar)^{-1} \sum_i (\feta_i dz_i - \fz_{i} d \eta_{i}) 
$$
gives a Fedosov connection, and flat sections of the Weyl bundle are given by
$$
F(z+\fz, \eta + \widehat {\eta})
$$
where $f$ is a holomorphic
function in a small polydisc. But this is precisely the Weyl deformation of
$\C^{2n}$. 

\subsection{Construction of Fedosov connections}

Let $(M,\omega )$ be a complex manifold with a holomorphic symplectic structure such that the
maps 
\begin{equation}
  \label{complex21}
 H^i (M ,\C )\rightarrow H^i(M,{\cal O}_M) ;\ i=1,2
\end{equation}
are surjective. Fix a splitting

\begin{equation}  \label{splitting of H2}
H^2(M) \stackrel{\sim}{\rightarrow} H^2(F^1 \Omega^{*,*} (M)) \oplus H^{0,2}(M)
\end{equation}
where $F^p \Omega^{*,*} (M) = \Omega^{\geq p,*}$
is the Hodge filtration.

\begin{th}\label{complex11}
Under the assumption above, let $\alpha \in \frac{1}{i\hbar}\omega + H^2(F^1 \Omega^{*,*} (M))[[\hbar]]$. There exists unique element $\tau \in \hbar H^{0,2}(M)[[\hbar]]$ such that $\alpha + \tau$ is a characteristic class of a Fedosov connection.
\end{th}
{\em Proof}.
\begin{lemma} \label{dbar+A-1}
The embedding 
$$ (\Omega^{*,*}(M, \C[[\hbar]]), \overline{\partial}) \rightarrow (\Omega^{*,*}(M, \tilde{\frak{g}}), \overline{\partial}+\mbox{ad}A_{-1})$$
is a quasi-isomorphism. In particular, the subcomplex $\sum_{p+q \; {odd}}(\Omega^{p,*}(M, \tilde{\frak{g}}_q))$ is acyclic with respect to $\overline{\partial}+\mbox{ad}A_{-1}$.
\end{lemma}

{\em{Proof}}. The lemma is implied by the fact that the differential $\mbox{ad}A_{-1}$ is acyclic in positive degrees, and its cohomology in degree zero is $\Omega^{0,*}(M)[[\hbar]]$.

  Now suppose we are given an element $(i\hbar)^{-1} \omega +\theta +\phi$
of $ ((i\hbar)^{-1} \omega + \Omega^{20}(M,\k ))\oplus \Omega^{11} (M,\k )$
representing a class in $ H^2 (F^1 \Omega^2 (M,\k ))$.

The construction of the Fedosov connection procedes by induction over the grading of the
Lie algebra $\tilde{\frak g}$ just as in the smooth case.

{\bf 1}. The flatness of $\nabla$ gives the following equations for
the pair $(A_0 , B_1)$:

\begin{eqnarray}
    \label{complex12}
  \nabla_0 A_{-1} +[A_{-1},A_0]=0 \\
  \dbar A_0 +[A_{-1},B_1 ]= \phi_0 \\
   \label{complex12'} \dbar B_1 = 0
\end{eqnarray}

Because of lemma \ref{dbar+A-1}, this system of equations has a solution $(A_0, B_1)$.

{\bf 2}.  
Given
($A_0 ,B_1 $), we want to find ($A_1, B_2 $) satisfying the equations
\begin{eqnarray}
  \label{complex14}
  [A_{-1},A_1] +\frac{1}{2} [\nabla_0 +A_0 ,\nabla_0 +A_0]=\theta_0 ,\\
  \mbox{$[ \dbar ,A_1 ]$} +[\nabla_0 +A_0 ,B_1 ] +[A_{-1},B_2 ]=0 ,\\
    \label{complex14'}
\dbar B_2 + \frac{1}{2}[B_1,B_1] = \tau_2       
\end{eqnarray}
The Bianchi identity for $\nabla_0 + ad A_{-1} +adA_{0}$
implies that $(\frac{1}{2} [\nabla_0 +A_0 ,\nabla_0 +A_0]-\theta_0 ,\;[\nabla_0 +A_0 ,B_1 ],\;\frac{1}{2}[B_1,B_1])$ is a $\dbar + \mbox{ad}A_{-1}$-cocycle. By lemma \ref{dbar+A-1}, the Dolbeault cohomology class of $\tau_2$ for which (\ref
{complex14}-\ref{complex14'}) has a solution $(A_1, \;B_2)$ exists and is unique.

Now, assume that we have already constructed the terms $A_{\leq 2n-1}$, $B_{\leq{2n}}$, and $\tau_{\leq 2n}$. One has
 \begin{eqnarray} \label{complex15}
  [A_{-1},A_{2n}] + [\nabla_0 +A_0 , A_{2n-1}] + \ldots  = 0,\\
  \mbox{$[ \dbar ,A_{2n} ]$} +[A_{-1}, B_{2n+1} ] +[\nabla_0 +A_0 ,B_{2n} ]+ \ldots = \phi_{2n} ,\\
   \dbar B_{2n+1} + [B_1, B_{2n}] + \ldots = 0      
\end{eqnarray}
By Bianchi identity for the connection $\nabla_0 + A_{\leq 2n-1} + B_{\leq 2n}$, and by lemma \ref{dbar+A-1}, this system has a solution $(A_{2n}, B_{2n+1})$. Now one has 
 \begin{eqnarray} \label{complex16}
  [A_{-1},A_{2n+1}] + [\nabla_0 +A_0 , A_{2n}] + \ldots  = \theta_{2n},\\
  \mbox{$[ \dbar ,A_{2n+1} ]$} +[A_{-1}, B_{2n+2} ] +[\nabla_0 +A_0 ,B_{2n+1} ]+ \ldots = 0 ,\\
   \label{complex16'}
 \dbar B_{2n+1} + [B_1, B_{2n}] + \ldots = \tau_{2n+2}      
\end{eqnarray}
By Bianchi identity for the connection $\nabla_0 + A_{\leq 2n} + B_{\leq 2n+1}$, and by lemma \ref{dbar+A-1}, there exists unique Dolbeault class $\tau_{2n+2}$ for which the system (\ref{complex16}-\ref{complex16'}) has a solution $(A_{2n}, B_{2n+1})$.

It remains to show that the cohomology classes $\tau_{2n+2}$ depend only on the cohomology classes of $\theta _{\leq 2n-2} $ and $\phi_{2n}$ (and not on the choices of $A_i, \; B_i$). But this follows immediately from remark \ref{gauge partial}.

\begin{th}
Let $(M,\omega )$ be a complex symplectic manifold such that the
maps 
\begin{equation}
  \label{complex21}
 H^i (M ,\C )\rightarrow H^i(M,{\cal O}_M) ;\ i=1,2
\end{equation}
are surjective. The set of isomorphism classes of formal deformations of $(M,\omega
)$ is in bijective correspondence with $H^2 (F^1 \Omega^{*,*}
(M),d)[[\hbar ]]$. Moreover there exists a family of smooth (nonlinear) maps:
\begin{equation}
  \label{cc10}
\tau_n : \{ H^2 (F^1 \Omega^{*,*} (M),d) \}^n \rightarrow H^2(M,{\cal O}_M) )
\end{equation}
such that the characteristic class of the deformation associated to $\alpha_0
+\hbar \alpha_1 +\ldots $ is given by the sum
\begin{equation}
  \label{eq:intcomp 2}
  \frac{1}{i\hbar}+\sum_n \hbar^n (\alpha_n + \tau_n (\alpha_0 , \ldots , \alpha_{n-1} )).
\end{equation}
The associated formal deformation of the sheaf of algebras of holomorphic
functions is locally isomorphic to the Weyl deformation of holomorphic
functions on an open polydisc in $\C^{dim M}$.
\end{th}
{\em Proof}. Follows immediately from combining theorems \ref{structure of formal deformations}, \ref{complex11}, and \ref{ccc5}.

\section{Index theorems}   \label{Index theorems}
\subsection{The trace density}   \label{The trace density}
Let $\fA _M$ be a deformation of a symplectic Lie algebroid $E$ with characteristic class $\theta$. Recall that for any unital algebra $A$ over $\C$ the Hochschild complex $(C_*(A,A), b)$ is defined, along with the negative cyclic complex $CC^{-}_*(A) = (C_*(A,A)[[u]], b+uB)$ and the periodic cyclic complex $CC^{per}_*(A) = (C_*(A,A)[ u^{-1},u]], b+uB)$ where $u$ is a formal parameter of degree $-2$; $C_n(A,A) = A\otimes (A/\C1)^{\otimes n}$; $b: C_{n}(A,A) \rightarrow  C_{n-1}(A,A)$ is the Hochschild differential, and $B: C_{n}(A,A) \rightarrow  C_{n+1}(A,A)$ is the cyclic differential (cf. \cite{Loday}).

Let $n = \frac{1}{2} \mbox{dim}E$. In this subsection, we construct {\em the trace density maps} 
\begin{eqnarray}
\label{trace C}
\mu^{\hbar}: C_*(\fA (M)) \rightarrow ({}^E\Omega ^{2n-*}(M)((\hbar)), 0)\\
\label{trace CC-}
\mu^{\hbar}: CC^-_*(\fA (M)) \rightarrow ({}^E\Omega ^{2n-*}(M)((\hbar))[[u]], d)\\
\label{trace CCper}
\mu^{\hbar}: CC_*^{per}(\fA (M)) \rightarrow ({}^E\Omega ^{2n-*}(M)((\hbar))[u^{-1},u]], d)
\end{eqnarray}
The construction is as follows. Let $\hat{\Omega}^*$ be the completion of the space of differential forms on ${\mathbb{R}}^n$ at $0$. This is a module over ${\frak{g}}=\mbox{Der} \vfA$ whose action is induced by the homomorphism of reduction modulo $\hbar$
$$  {\frak{g}} \rightarrow \mbox{Ham}({\mathbb{R}}^n)$$
where $\mbox{Ham}({\mathbb{R}}^n)$ is the algebra of formal Hamiltonian vector fields.
Therefore
\begin{equation}  \label{module of homs}
{\mathbb{L}}^* = \mbox{Hom}(CC^-_*(\vfA), \hat{\Omega}^{2n-*})
\end{equation}
is a complex of ${\frak{g}}$-modules.

In \cite{BNT}, we constructed the canonical element $\mu^{\hbar}$ of degree zero in ${\mathbb{L}}^*$. The image of $\mu^{\hbar}$ under the Gelfand-Fuks map (definition \ref{gf}) induces a map of sheaves
\begin{equation}  \label{tr on sheaves}
{}^E\Omega^*_M(CC^-_*(\vfA)) \rightarrow {}^E\Omega^*_M(\hat{\Omega}^{2n-*})((\hbar))
\end{equation} 
Since the complexes of sheaves ${}^E\Omega^*_M(CC^-_*(\vfA))$ and $CC^-_*(\fA _M)$, resp. ${}^E\Omega^*_M(\hat{\Omega}^{2n-*})((\hbar))$ and ${\Omega}_M^{2n-*}((\hbar))$, are quasi-isomorphic,
one gets the map (\ref{trace CC-}). To get the map (\ref{trace C}), one puts $u=0$, and to get the map (\ref{trace CCper}), one localizes with respect to $u$. Cf. \cite{BNT} for details, including the explicit definition of the module ${\mathbb{L}}^* $.

\subsection{Index theorem for symplectic Lie algebroids} \label{index for E}
Define the $\C[[u]]$-linear continuous morphism 
$$\mu: \;  CC_*^{per}(\fA _M) \rightarrow ({}^E \Omega^*_M [[u]], ud)$$
as the projection 
$$CC_*^{per}(\fA _M) \rightarrow CC_*^{per}(C^{\infty}_M)$$
followed by the Connes' quasi-isomorphism
$$a_0 \otimes \ldots \otimes a_p \mapsto \frac{1}{p!}a_0 da_1 \ldots d a_p$$

Let
$$j: ({}^E\Omega ^{2n-*}_M[[u]],d) \rightarrow ({}^E\Omega ^{*}_M[[u]],ud)$$
be the morphism of complexes given by 
$$j(u^p \alpha) = u^{p+k-n}\alpha$$
for $\alpha \in {}^E\Omega ^{p}_M$. Let $\hat{A}(E)$ be the $\hat{A}$ class of a $U(n)$-bundle obtained by reducing the structure group to its maximal compact subgroup. Finally, by
$$i: \Omega_M^* \rightarrow {}^E\Omega^*_M ((\hbar))$$
we denote the composition of the map $i: \Omega_M^* \rightarrow {}^E\Omega^*_M $
conjugate to the anchor map $\rho$ with the embedding ${}^E \Omega_M^* \hookrightarrow {}^E\Omega^*_M((\hbar)) $.
\begin{th}   \label{index thm for E}
$$j\circ\mu^{\hbar}=\sum_{p\geq 0}u^p (\hat{A}(E)e^{\theta})_{2p}(i\circ \mu)$$
\end{th}
{\em Proof}. Follows immediately from the Riemann-Roch theorem for periodic cyclic cochains \cite{BNT}.


\bigskip
$^1$ Mathematics Institute 

University of Copenhagen

Universitatsparken 5 2100 Copenhagen, Denmark

rnest@math.ku.dk

\bigskip

$^2$ Department of Mathematics

Penn State University

University Park PA 16802 USA

tsygan@math.psu.edu

Supported
in part by NSF grant DMS-9504522.

\end{document}